\documentclass[preprint,5p]{elsarticle}
\usepackage{multicol}
\usepackage[utf8]{inputenc}
\usepackage[english]{babel}
\usepackage{amsthm,amsmath,amssymb,mathtools}
\usepackage{graphicx}
\usepackage{bbm}
\usepackage{floatrow}
\usepackage{hyperref}
\usepackage{cleveref}
\usepackage{lineno}
\usepackage{caption}
\usepackage{subcaption}

\usepackage{float}

\newtheorem{theorem}{Theorem}
\newtheorem*{theorem*}{Theorem}
\newtheorem{lemma}[theorem]{Lemma}
\newtheorem{corollary}[theorem]{Corollary}
\newtheorem{proposition}[theorem]{Proposition}
\newtheorem{ass}{Assumption}

\theoremstyle{remark}
\newtheorem{remark}{Remark}
\newtheorem{deff}[remark]{Definition}
\newtheorem*{probdef*}{Problem-Definition}
\newtheorem*{deff*}{Definition}

\numberwithin{equation}{section}

\newcommand{\bfx}{\mathbf{x}}

\begin{document}

\begin{frontmatter}

\title{Interplay between depth and width for interpolation in neural ODEs}
\author[1]{Antonio Álvarez-López\corref{cor1}}
\ead{antonio.alvarezl@uam.es}
\author[2]{Arselane Hadj Slimane}
\ead{arselane.hadj_slimane@ens-paris-saclay.fr}
\author[1,3,4]{Enrique Zuazua}
\ead{enrique.zuazua@fau.de}
\cortext[cor1]{Corresponding author}
\affiliation[1]{organization={Universidad Autónoma de Madrid},
addressline={C. Francisco Tomás y Valiente, 7},
city={Madrid},
postcode={28049},
country={Spain}}
\affiliation[2]{organization={ENS Paris Saclay},
addressline={Avenue des science, 4},
city={Gif-sur-Yvette},
postcode={91190},
country={France}}
\affiliation[3]{organization={Chair for Dynamics, Control, Machine Learning, and Numerics, Alexander von Humboldt-Professorship, Department of Mathematics, Friedrich-Alexander-Universität Erlangen-Nürnberg},
addressline={Cauerstraße, 11},
city={Erlangen},
postcode={91058},
country={Germany}}
\affiliation[4]{organization={Fundación Deusto},
addressline={Av. de Universidades, 24},
city={Bilbao},
postcode={48007},
country={Spain}}

\begin{abstract}
Neural ordinary differential equations (neural ODEs) have emerged as a natural tool for supervised learning from a control perspective, yet a complete understanding of their optimal architecture remains elusive. In this work, we examine the interplay between their width $p$ and number of layer transitions $L$ (effectively the depth $L+1$).  Specifically, we assess the model expressivity in terms of its capacity to interpolate either a finite dataset $\mathcal{D}$ comprising $N$ pairs of points or two probability measures in $\mathbb{R}^d$ within a Wasserstein error margin $\varepsilon>0$. Our findings reveal a balancing trade-off between $p$ and $L$, with $L$ scaling as $O(1+N/p)$ for dataset interpolation, and $L=O\left(1+(p\varepsilon^d)^{-1}\right)$ for measure interpolation. 

In the autonomous case, where $L=0$, a separate study is required, which we undertake focusing on dataset interpolation. We address the relaxed problem of $\varepsilon$-approximate controllability and establish an error decay of $\varepsilon\sim O(\log(p)p^{-1/d})$. This decay rate is a consequence of applying a universal approximation theorem to a custom-built Lipschitz vector field that interpolates $\mathcal{D}$. In the high-dimensional setting, we further demonstrate that $p=O(N)$ neurons are likely sufficient to achieve exact control.
\end{abstract}
\begin{keyword}
Neural ODEs  \sep Depth \sep Width \sep Simultaneous controllability \sep Transport equation \sep Wasserstein distance

\MSC 34H05 \sep 35Q49 \sep 68T07 \sep 93B05 
\end{keyword}
\end{frontmatter}

\section{Introduction}

Residual neural networks (ResNets) are formally defined as the family of discrete systems
\begin{align}\label{eq:resnet-intro}
\begin{cases}
\mathbf{x}_{k+1}&=\mathbf{x}_k+W_k\boldsymbol{\sigma}\left( A_k \mathbf{x}_k+\mathbf{b}_k\right), \\
\mathbf{x}_0&\in \mathbb{R}^d,
\end{cases}
\end{align}
 where $k =0, \ldots, L$,
$W_k\in\mathbb{R}^{d\times p}$, $A_k\in\mathbb{R}^{p\times d}$ and $\mathbf{b}_k\in\mathbb{R}^p$, for some $d\geq~1$, $L\geq0$ and $p\geq1$. Each time step $k$ identifies a layer of the network.
The number of layers $L+1$ is the \emph{depth} of \eqref{eq:resnet-intro}. The parameter $p$ is the \emph{width} of \eqref{eq:resnet-intro}, identifying the number of neurons per layer. The activation function $\boldsymbol{\sigma}: \mathbb{R}^p \to \mathbb{R}^p$ is defined as the column vector $\boldsymbol{\sigma}(\mathbf{y})=\big(\sigma(y^{(1)}),\dots,\sigma(y^{(p)})\big)^\top$ from a chosen nonlinear function $\sigma:~\mathbb{R}\to\mathbb{R}$. We consider the Rectified Linear Unit (ReLU), given by
$\sigma(z)=\max\{z,0\},\;z\in \mathbb{R}.$

It has been noted \cite{weinan,haber,chen2019neural,chang2018multilevel,Ruthotto2024differential} 
that \eqref{eq:resnet-intro} can be identified with the forward Euler discretization scheme for the class of continuous models known as neural ordinary differential equations (neural ODEs),
\begin{align} \label{eq:node-intro}
\begin{cases}
\dot\bfx(t)&= W(t) \boldsymbol{\sigma} \left(A(t)\, \bfx(t)+\mathbf{b}(t)\right),\\
\bfx(0)&=\bfx_0 \in \mathbb{R}^d,
\end{cases}
\end{align}
where $(W,A,\mathbf{b})\in L^\infty\left((0,T),\mathbb{R}^{d\times p}\times\mathbb{R}^{p\times d}\times\mathbb{R}^p\right)$ for some $T>0$.
Here, $t \in (0, T)$ parameterizes the evolution of the states through a continuous range of layers. As discussed in \cite{rbz-node,alvlop24}, it is common to assume that $(W, A, \mathbf{b})$ is a step function over $(0, T)$, to align closer with the dynamics of \eqref{eq:resnet-intro}. Then, since $\sigma$ is Lipschitz, existence and uniqueness of solutions hold for any $(W, A, \mathbf{b})$ and initial condition $\mathbf{x}_0$. Equation \eqref{eq:node-intro} can be equivalently written as
\begin{equation}\label{eq:node-p}
\dot\bfx=\sum_{i=1}^p \mathbf{w}_i(t)\sigma(\mathbf{a}_i(t)\cdot\bfx+b_i(t)),
 \end{equation}
where $\mathbf{w}_i$ and $\mathbf{a}_i$ are respectively the $p$ columns of $W$ and the $p$ rows of $A$, both seen as column vectors in $\mathbb{R}^d$, while $b_i$ is the $i$-th coordinate of $\mathbf{b}$, for $i=1,\dots,p$. In this work, we use formulation \eqref{eq:node-p}, although, for simplicity, we represent $(\mathbf{w}_i,\mathbf{a}_i,b_i)_{i=1}^p$ in their matrix form $(W,A,\mathbf{b})$, which corresponds to the equivalent system \eqref{eq:node-intro}.
\begin{figure}[t]  
\centering
\begin{subfigure}{0.49\linewidth}
    \includegraphics[width=\linewidth]{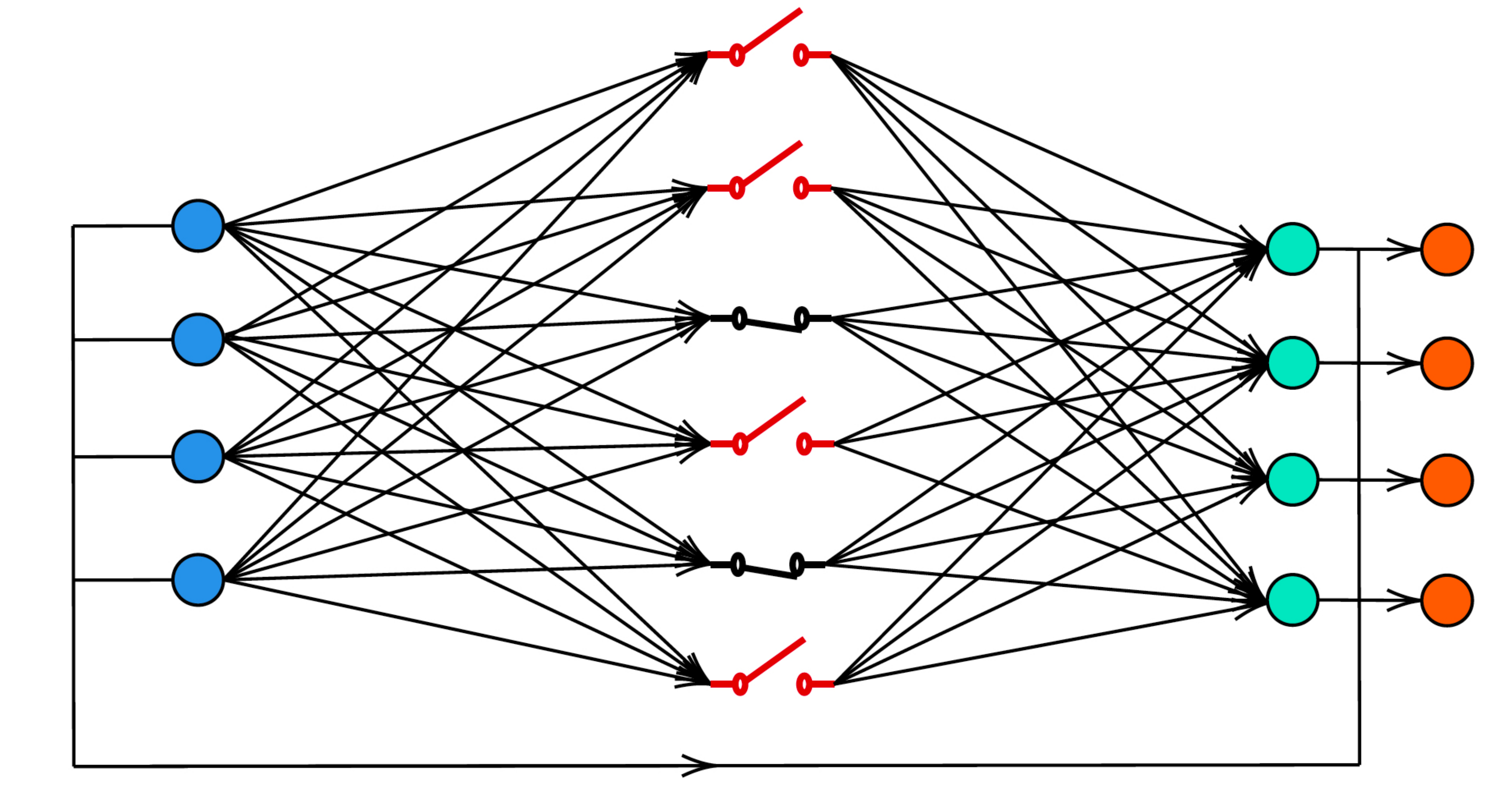}
   \caption{Shallow ResNet}
    \label{fig:sresnet} 
\end{subfigure}
\begin{subfigure}{0.49\linewidth}
    \includegraphics[width=\linewidth]{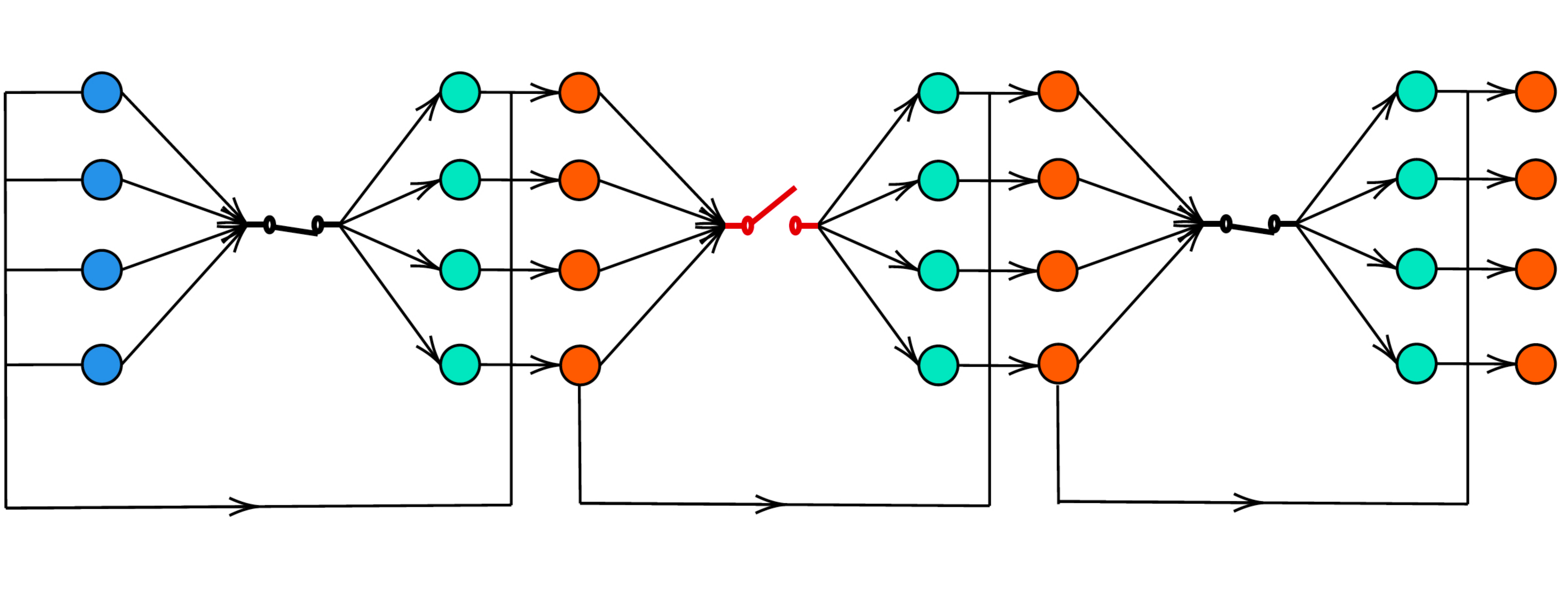}
     \caption{Narrow ResNet}
    \label{fig:tresnet} 
\end{subfigure}
\caption{Qualitative representation of models \eqref{eq:node-shallow} and \eqref{eq:node-Narrow} as discrete systems. Blue circles represent the input $\mathbf{x}$; switches depict ReLU functions; green circles indicate the result of  $W\boldsymbol{\sigma}(A\mathbf{x}+\mathbf{b})$; orange circles represent the output after residual connections.}
\end{figure}
Equation \eqref{eq:node-p} can be naturally extended to handle probability distributions, rather than points in $\mathbb{R}^d$, by interpreting its right-hand side as the advection field that drives the evolution of a measure $\rho$. This extension gives rise to the neural transport equation \cite{domzuazuaNormFlows}:
\begin{equation}\label{eq:neuraltransport}
    \partial_t \rho +\operatorname{div}_\bfx\Big(\sum_{i=1}^p \mathbf{w}_i\,\sigma(\mathbf{a}_i\cdot\bfx+b_i)\rho\Big)=0.
\end{equation} 
Prior research indicates that control theory offers significant potential for examining the properties of neural ODEs, for instance, via optimal control \cite{weinan,largetime,Yage2021SparseAI} or geometric control techniques \cite{tabuada23,Agrachev2021ControlOT,scagliottidiffeo}. A fundamental problem still open is to develop a comprehensive understanding of the roles played by depth and width with respect to the expressive power of the model, see \cite{hardt2017identity,expressivepower17,quasiequiv}. This property is often evaluated by its capacity to interpolate either a finite set of point pairs or two given probability measures. 

The first, commonly referred to as finite-sample expressivity \cite{Yun2018SmallRN}, is associated with the approximation power of the model, see \cite{qianxiao22}. It essentially amounts to a simultaneous control problem, where the aim is to find a control function $(W, A, \mathbf{b})$ such that the associated input-output map, given by the flow of \eqref{eq:node-p}, maps $N$ specified points to their $N$ corresponding target points in $\mathbb{R}^d$. Throughout this work, figures represent each input point as a colored solid circle and its corresponding target point as an empty circle of the same color.

In the second scenario, we aim to control the transport dynamics described by \eqref{eq:neuraltransport} in order to transform a given initial density $\rho_0$ into another density $\rho_T$. This task is highly relevant for probabilistic modeling or the generation of synthetic data via normalizing flows \cite{kobyzev,normflowspapamakarios,Grathwohl}. We approach it as an approximate control problem in the Wasserstein-$q$ metric space for $q\geq1$, extending prior work focused on $W_2$, see \cite{osher}, or in $W_1$ with $p=1$, see \cite{rbz-node}.

Our main objective is the development of a comprehensive theory of interpolation for the family of models described by \eqref{eq:node-p}, linking the error to the specific architecture given by $p$ and $L$. Both numerical and theoretical studies \cite{huang2003,powerofdepth,mhaskar2017and} suggest that networks with greater depth usually achieve better performance. This tendency is particularly noticeable in training \cite{understandingdl,Yun2018SmallRN}, which reinforces the intuition that a deeper network should possess greater expressivity, i.e., an enhanced ability to learn more complex non-linear functions. Understanding the balance between width and depth is thus vital for the optimal design of networks. We tackle this significant question using the continuous framework of neural ODEs, where depth is expressed as $L+1$,  $L$ being the number of time discontinuities of the control $(W,A,\mathbf{b})$. As we vary $L$ and $p$, two limiting models emerge.\newline
\textbf{Shallow neural ODEs.} $L=0$ is fixed, while the width $p$ can be as large as required:
    \begin{equation}\label{eq:node-shallow}    
        \dot\bfx=\sum_{i=1}^p \mathbf{w}_i\sigma\left( \mathbf{a}_i\cdot\bfx+b_i\right),
    \end{equation}
    where $\{(\mathbf{w}_i,\mathbf{a}_i,b_i)\}_{i=1}^p\subset\mathbb{R}^d\times\mathbb{R}^d\times\mathbb{R}$ are constant controls, making the equation autonomous. The field on the right-hand side of \eqref{eq:node-shallow} corresponds to a one hidden layer neural network with $d$ components. The approximation capacity of this class of functions has been extensively studied (see \cite{cybenko,pinkus_1999,devore_hanin_petrova_2021}). The discrete version of \eqref{eq:node-shallow} can be identified with a one hidden layer ResNet (see \cref{fig:sresnet}).\newline
\textbf{Narrow neural ODEs.} $p=1$ is fixed, while the depth $L+1$ can be as large as required:
    \begin{equation}\label{eq:node-Narrow}
        \dot\bfx= \mathbf{w}(t)\sigma\left(\mathbf{a}(t)\cdot\bfx+b(t)\right),
    \end{equation}
    where $(\mathbf{w},\mathbf{a},b)\in L^\infty\left((0,T),\mathbb{R}^d\times\mathbb{R}^d\times\mathbb{R}\right)$. The ability of this model to interpolate data and approximate functions has been explored in \cite{rbz-node,uat1neuronresnet}. It offers the advantage of easier construction of explicit controls compared to \eqref{eq:node-shallow}, owing to its simplified dynamics, albeit at the expense of increased depth, which scales with the cardinal $N$ of the dataset. The discrete version of \eqref{eq:node-Narrow} corresponds to a deep ResNet with one neuron per hidden layer, so it alternates layers of dimension 1 and $d$ (see \cref{fig:sresnet}).

Developing a unified theory that bridges shallow and narrow neural ODEs would combine the vast work done for \eqref{eq:node-shallow} with the intuitive dynamics of \eqref{eq:node-Narrow}. Moreover, it would facilitate the optimal design of a neural ODE through the strategic choice of depth and width. This entails optimizing the complexity $\kappa$, defined as the total number of parameters in \eqref{eq:node-p}:
\begin{equation}
    \label{eq:complexity}
    \kappa\coloneqq(L+1)\times p\times(2d+1).
\end{equation}
Indeed, on each of the $L+1$ hidden layers, $p$ neurons of dimension $2d+1$ need to be determined. 

\subsection{Roadmap}
In \cref{sec:mainres}, we present the main results of our work in two parts. First, in \cref{subsec:simcontrol}, we study the problem of interpolating a finite dataset in \eqref{eq:node-p}, which is recast as the property of simultaneous control for neural ODEs.  Second, in \cref{subsec:transp}, we approach approximate controllability of probability measures using the dynamics provided by the neural transport equation \eqref{eq:neuraltransport}. In \cref{conclusions}, we discuss the main implications of our work and pose some open questions. In \cref{sec:proofs}, we prove the main results and provide the necessary tools as lemmas.

\subsection{Notation}\label{subsec:notations}
\begin{itemize}
\item We use subscripts to identify the particular elements from a dataset and superscripts for the coordinates of a vector. In addition, (column) vectors are denoted with bold letters and matrices with capital letters.
    \item We denote by $\bfx\cdot\mathbf{y}$ the scalar product of $\bfx,\mathbf{y}\in\mathbb{R}^d$.
    \item We denote by $\lceil z\rceil$ the lowest integer greater than or equal to $z\in\mathbb{R}$, and by $\lfloor z\rfloor$ the highest integer lower than or equal to $z$.
    \item We denote by $\mathbb{S}^{d-1}$ the $(d-1)$-dimensional sphere in $\mathbb{R}^d$.
    \item We denote by $\operatorname{Lip}\left(\mathbb{R}^d,\mathbb{R}^d\right)$ the space of Lipschitz-continuous vector fields in the usual norm, and  by $L_V$ the Lipschitz constant of each  $\mathbf{V}\in\operatorname{Lip}\left(\mathbb{R}^d,\mathbb{R}^d\right)$.
    \item Given any Borel measure $\mu$ in $\mathbb{R}^d$ and any measurable function $f:\mathbb{R}^d\rightarrow\mathbb{R}^d$, we denote by $f_\#\mu$ the pushforward measure, defined for every Borel subset $A\subset\mathbb{R}^d$ by \begin{equation*}f_\#\mu(A)=\mu(f^{-1}(A)).\end{equation*}
\end{itemize}

\section{Main results}\label{sec:mainres}
\subsection{Simultaneous control}\label{subsec:simcontrol}
Let $N\geq1$, $d\geq 2$, and consider a dataset
\begin{equation}\label{sample}
\mathcal{D}=\{(\bfx_n,\mathbf{y}_n)\}_{n=1}^N\subset\mathbb{R}^d\times\mathbb{R}^d
\end{equation}
with $\bfx_n\neq\bfx_m$ and $\mathbf{y}_n\neq\mathbf{y}_m$ for all $n\neq m$. First, we study the finite-sample expressivity of the general model \eqref{eq:node-p}, recast as a problem of simultaneous control.
 \begin{probdef*}
     For any fixed time horizon $T>0$, find controls \begin{equation*}\left\{(\mathbf{w}_i,\mathbf{a}_i,b_i)\right\}_{i=1}^p\subset L^\infty\left((0,T);\mathbb{R}^d\times\mathbb{R}^d\times\mathbb{R}\right),\end{equation*}for some $p\geq1$, such that the flow $\Phi_T(\cdot;W,A,\mathbf{b})$ of \eqref{eq:node-intro} \emph{interpolates} the dataset $\mathcal{D}$, i.e., it simultaneously drives each data point from its initial position $\bfx_n$ to its target $\mathbf{y}_n$. This is fulfilled when
\begin{equation*}
\Phi_T(\bfx_n;W,A,\mathbf{b})=\mathbf{y}_n\qquad \text{ for all }n=1,\dots,N,
\end{equation*}
where $W,A,\mathbf{b}$ are respectively the matrix with columns $\mathbf{w}_i$, the matrix with rows $\mathbf{a}_i$ and the vector with components $b_i$, for $i=1,\dots,p.$
 \end{probdef*} 
Our first result provides a relationship between $L$ and $p$ that ensures interpolation of $\mathcal{D}$:
        \begin{theorem}\label{th:SimControl-p}
           Let $N\geq1$, $d\geq2$ and $T>0$ be fixed. Consider the dataset $\mathcal{D}$ as defined in \eqref{sample}.     For any $p\geq1$, there exists a piecewise constant control \begin{equation*}\left(W,A,\mathbf{b}\right)\in L^\infty\left((0,T);\mathbb{R}^{p\times d}\times\mathbb{R}^{p\times d}\times\mathbb{R}^p\right)\end{equation*}
           such that the flow $\Phi_T(\cdot;W,A,\mathbf{b})$ generated by \eqref{eq:node-p} interpolates the dataset $\mathcal{D}$, i.e., \begin{equation*}\Phi_T(\bfx_n;W,A,\mathbf{b})=\mathbf{y}_n,\qquad\text{for all }n=1,\dots,N.\end{equation*}
             Furthermore, the number of discontinuities of $(W,A,\mathbf{b})$ is
\begin{equation}\label{eq:Lvsp}L=2\left\lceil N/p \right\rceil-1.\end{equation}
        \end{theorem}

      \begin{remark} 
            If the target points $\{\mathbf{y}_n\}_{n=1}^N$ in \eqref{sample} are not distinct, interpolation is not achievable due to the uniqueness of solutions in the system \eqref{eq:node-p}. In such cases, we relax the statement from exact to approximate controllability by applying \cref{th:SimControl-p} to an $\varepsilon$-perturbation of the targets, for some $\varepsilon>0$.
            \end{remark}
            Let us briefly describe the algorithm. First, we pivot around the $x^{(1)}$-coordinate and control the remaining $d-1$ coordinates. Consequently, the trajectory of each data point $\bfx_n$ is confined within the hyperplane defined by the equation $x_{n}^{(1)}=x^{(1)}_n$. Then, we pivot using the controlled coordinates to adjust $x^{(1)}$. This algorithm requires a depth of $2\lceil N/p\rceil$ layers, which is independent of the dimension $d$, since a constant control suffices to simultaneously steer $d-1$ coordinates in the first step, assuming that $p>N$. \begin{remark}\label{rem:extthm1}
        Our approach is broadly applicable to any activation function, provided it meets the following three conditions:
\begin{equation*}
    \text{1. } \sigma \text{ loc. Lipschitz;}\quad \text{2. } \sigma(z)\vert_{z\leq0}=0;\quad \text{3. } \sigma(z)\vert_{z>0}>0.
\end{equation*}
This generalization guarantees the extension of \cref{th:SimControl-p} to more general activation functions such as the ReLU powers $\sigma^k(z)=\operatorname{max}\{z,0\}^k$, for $z\in\mathbb{R}$ and $k\geq1$, whose approximation properties have been recently studied in \cite{CABANILLA2024106073}.
            \end{remark}

   In \eqref{eq:Lvsp} we can see that, as the width $p$ increases, the number of discontinuities $L$ decreases with the same rate, meaning that width and depth play a similar role in the steering.   Nevertheless, a result on the optimal design of our interpolating models can be derived:
     \begin{corollary}\label{cor:optimalcomp}
       For the family of controls given by  \cref{th:SimControl-p} that ensure interpolation of $\mathcal{D}$, the minimal complexity is  \begin{equation*} \kappa_{\text{min}}=(4d+2)(N+1), \end{equation*} obtained when $p=1$, i.e., when the neural ODE belongs to the narrow model \eqref{eq:node-Narrow}. 
        \end{corollary}
The complete transition from the narrow model \eqref{eq:node-Narrow} to the shallow model \eqref{eq:node-shallow}, characterized by $L=0$ is not attained in \eqref{eq:Lvsp}.  Due to the division into two steps in the proposed algorithm, whenever $p>N$ the selected control will exhibit a single switch ($L=1$), reaching a two-layer architecture, rather than the autonomous ansatz \eqref{eq:node-shallow}. The restriction naturally raises the question of whether simultaneous control is possible in shallow neural ODEs \eqref{eq:node-shallow}. For this task, a reconsideration of the algorithm presented in \cite{rbz-node} becomes necessary.  In the high dimensional setting, and more precisely, when the dimension exceeds the number of data points ($d>N$), we can refine the statement of \cref{th:SimControl-p} to include the case $L=0$:    
\begin{corollary}\label{cor:change-basis}
           Let $N\geq1$, $d\geq2$ with $d>N$, and $T>0$ be fixed. Consider the dataset $\mathcal{D}$ as defined in \eqref{sample}.     For any $p\geq1$, there exists a piecewise constant control \begin{equation*}
               \left(W,A,\mathbf{b}\right)\in L^\infty\left((0,T);\mathbb{R}^{d\times p}\times\mathbb{R}^{p\times d}\times\mathbb{R}^p\right)
           \end{equation*} such that the flow $\Phi_T(\cdot;W,A,\mathbf{b})$ generated by \eqref{eq:node-p} interpolates the dataset $\mathcal{D}$.  Furthermore, the number of discontinuities of $(W,a,\mathbf{b})$ is
\begin{equation*}L=2\big(\left\lceil N/p \right\rceil-1\big).\end{equation*}
\end{corollary}
 The key idea is that, when $d>N$, the first step in the proof of \cref{th:SimControl-p} can be suppressed. This is done by transforming the $x^{(1)}-$axis so that each $\bfx_n$ shares the same first coordinate with $\mathbf{y}_n$, or, equivalently, by optimally reorienting the hyperplanes represented in  \cref{fig:node-control1}. For more insights on the proof, see \cref{fig:dgtrn} in \cref{sec:proofs}.

\begin{remark}
\Cref{cor:change-basis} suggests ideas similar to those in \cite[Theorem 5.1]{largetime}. In that result, interpolation is established for $d \geq N$ in a simplified neural ODE, when $\sup_{n=1,\dots,N} |\bfx_n - \mathbf{y}_n| < \varepsilon$ for a sufficiently small $\varepsilon > 0$, under a geometric assumption on the images of the targets through $\boldsymbol{\sigma}$. Moreover, an estimation of the control cost is obtained, which is linear with respect to $\varepsilon$. The generation of new synthetic coordinates until $d \geq N$ is not typically a problem, as discussed in \cite{augmentednodes}, where the technique of embedding the dataset in $\mathbb{R}^d \times \{0,\dots,0\}$ is proposed and its computational advantages are studied.
\end{remark}

\begin{figure}[t!]  
\begin{subfigure}{4cm}
\centering
    \includegraphics[scale=0.1]{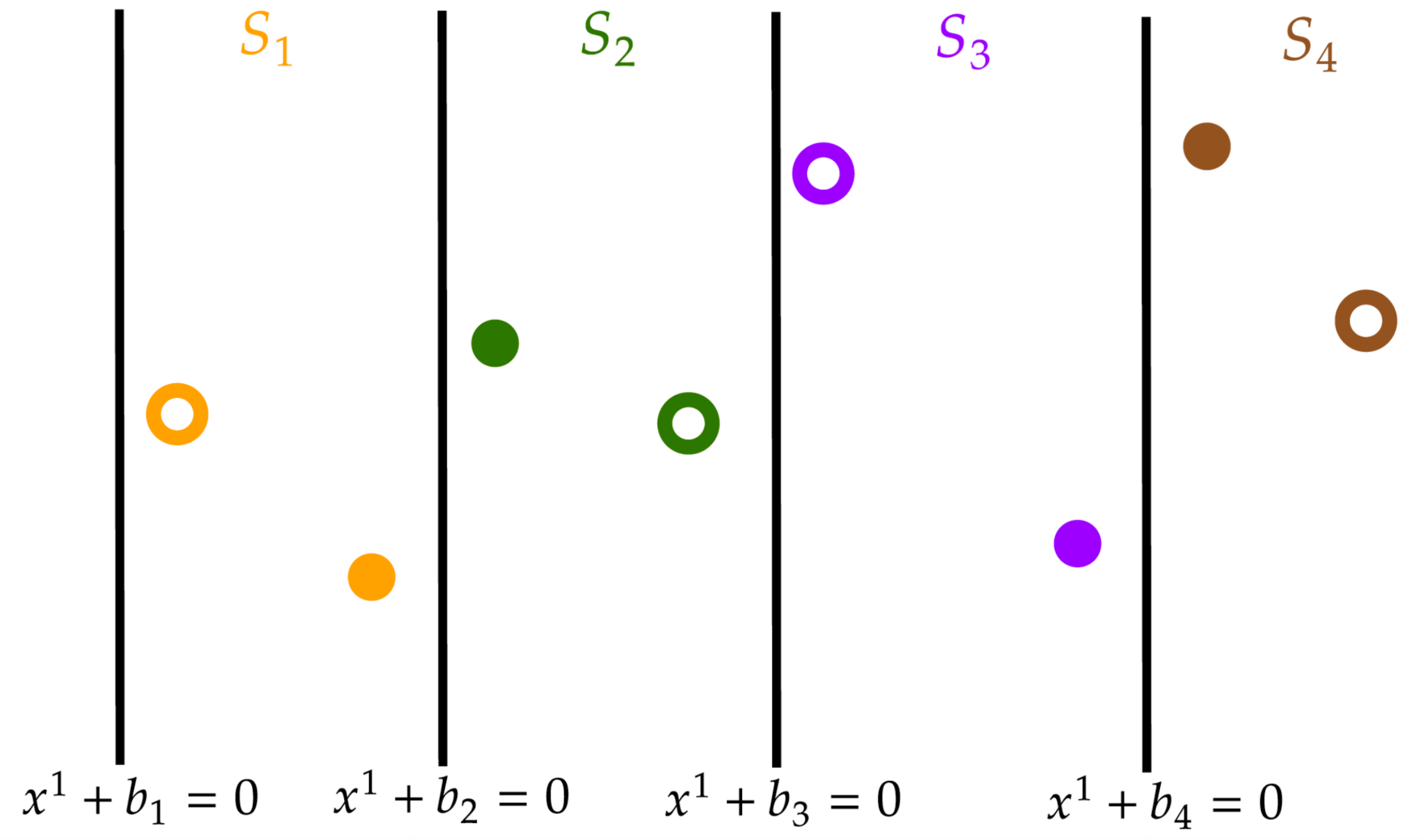}
    \label{fig:cybenko-sep}
\end{subfigure}
\;
\begin{subfigure}{4cm}
\centering
    \includegraphics[scale=0.1]{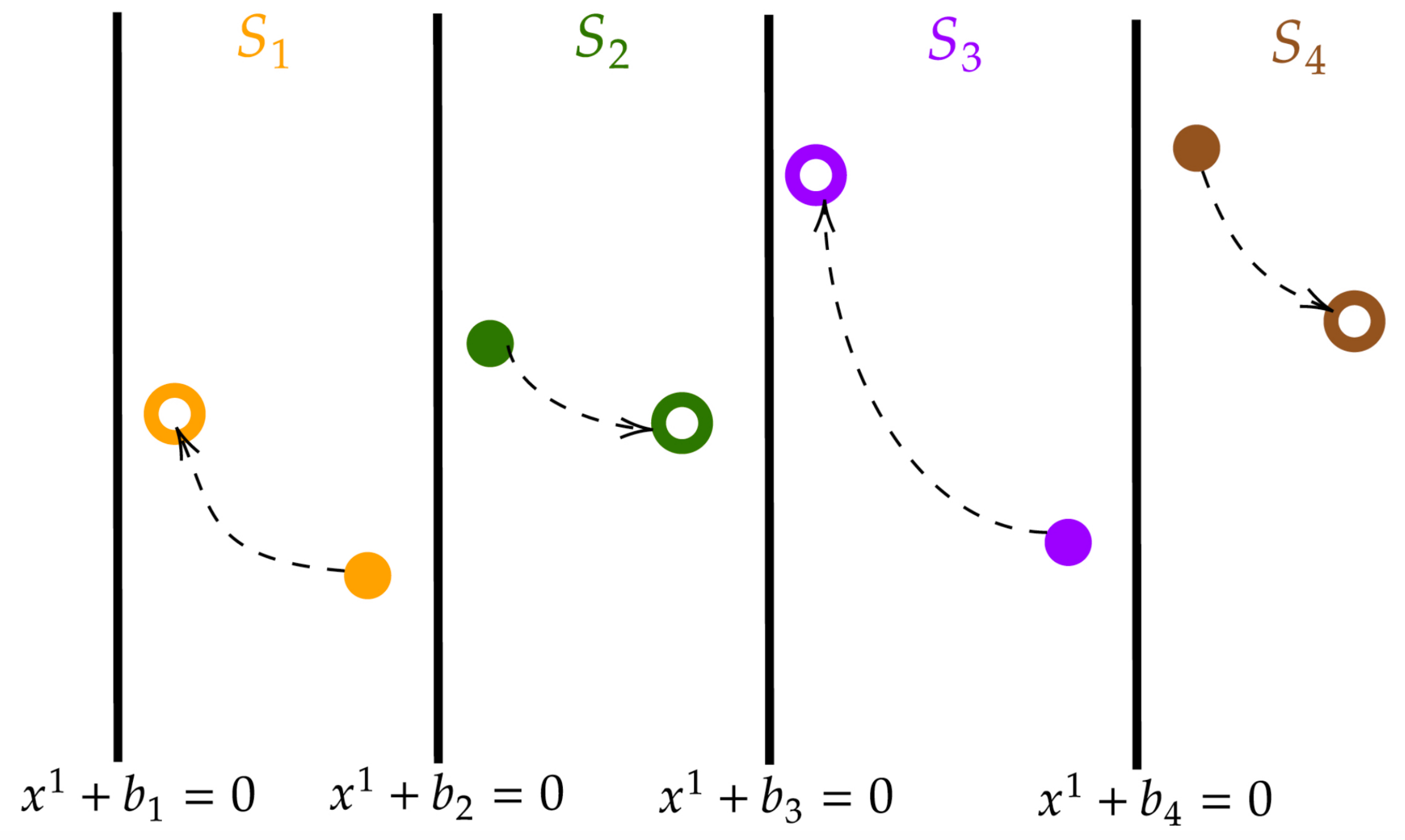}
    \label{fig:cybenko-control}
\end{subfigure}
\caption{Left: separability condition in \cref{hyp1}, for $\mathbf{a}=\mathbf{e}_1$. Right: trajectories for exact control in the same example.}
\end{figure}

In practice, $N$ tends to be larger than $d$. In that case, interpolation with constant controls can be obtained for $p=N$ under a certain separability hypothesis on $\mathcal{D}$:
        \begin{ass}\label{hyp1}
Let $\mathcal{D}=\{(\bfx_n,\mathbf{y}_n)\}_{n=1}^N\subset\mathbb{R}^d\times\mathbb{R}^d$ as defined in \eqref{sample}. There exist a vector $\mathbf{a}\in\mathbb{S}^{d-1}$, a permutation $\tau$ of $N$ elements and a sequence $-\infty<b_{N+1}<b_N<\cdots<b_{1}<\infty$ such that, for all $n=1,\dots,N-1$,
\begin{equation*}
-b_n<\mathbf{a}\cdot\bfx_{\tau(n)}<-b_{n+1}\quad \text{and}\quad-b_n<\mathbf{a}\cdot\mathbf{y}_{\tau(n)}<-b_{n+1}.\end{equation*}
\end{ass}
\Cref{hyp1} claims that we can diagonally separate each pair $(\bfx_n,\mathbf{y}_n)$ from the rest, in the sense that we can define $N+1$ parallel hyperplanes $H_n=\left\{\mathbf{a}\cdot\bfx+b_n=0\right\}$ such that the strip $S_n$ bounded by $H_n$ and $H_{n+1}$ contains only the point $\bfx_n$ and its target $\mathbf{y}_n$, for $n=1,\dots,N$ (see \cref{fig:cybenko-sep}). While the hypothesis might seem overly restrictive, it is noteworthy that if the points are randomly sampled from a compact set, the probability that the condition is fulfilled converges to 1 when the dimension grows:
 \begin{proposition}\label{lem:probability}
Let $\mu\in\mathcal{P}_{ac}^c(\mathbb{R}^d)$ such that the random variables $\pi_iX$ are independent and identically distributed (i.i.d.) for $i =1,\ldots,d$, where $X\sim\mu$ and $\pi_i$ is the canonical projection on the $i$-th coordinate. If every $\mathbf{x}_n$ and $\mathbf{y}_n$ in $\mathcal{D}$ is sampled from $\mu$, and $N$ is sufficiently large, then the probability $P$ that \cref{hyp1} is satisfied is bounded as
    \begin{equation*}
1-\left[1-\frac{1}{\sqrt{2}}\left(\frac{e}{2N}\right)^N\right]^d\leq P\leq 1.
\end{equation*}
        \end{proposition}
            The hypothesis of \cref{lem:probability} are fulfilled by the uniform probability measure in any hypercube, or by any isotropic Gaussian distribution.

Now, under \cref{hyp1}, we can build a constant control such that the flow of \eqref{eq:node-shallow}, taking a width $p=N$, interpolates $\mathcal{D}$. Our result is somehow a dynamic version of \cite[Theorem 1]{understandingdl}, under certain geometric conditions. A representation of the trajectories can be seen in \cref{fig:cybenko-control}. 
        \begin{corollary}
        \label{cor:snodewithsep}
           Consider a dataset $\mathcal{D}\subset\mathbb{R}^d\times\mathbb{R}^d$ for $d\geq2$, under \cref{hyp1}.  For any fixed $T>0$, there exists a control  $(W,A,\mathbf{b})\in \mathbb{R}^{d\times N}\times\mathbb{R}^{N\times d}\times\mathbb{R}^N$ such that the flow $\Phi_T$ generated by \eqref{eq:node-shallow} interpolates the dataset $\mathcal{D}$.
        \end{corollary}

All in all, new strategies are required to study simultaneous control in the autonomous model \eqref{eq:node-shallow} under general conditions. A natural starting point to assess the problem's feasibility is to relax it by admitting an error $\varepsilon > 0$, and use density tools provided by universal approximation theorems (UATs); see \cite{cybenko, pinkus_1999}. In deep learning, UATs establish the density of neural networks in function spaces over compact domains. The decay rate of the error in relation to the number of parameters of the network has been quantified for certain spaces \cite{bachbreaking, devore_hanin_petrova_2021}. Specifically, these studies bound the uniform error decay rate when the target function is Lipschitz continuous in a compact domain.

We will approach the UAT, often interpreted in a static manner, from our dynamic control perspective. In this regard, shallow neural ODEs provide a vector field that transitions initial data to final data, with its flow at time $T$ approximating the target function. First, we establish the existence of a time-independent field whose integral curves guide each input point $\bfx_n$ in $\mathcal{D}$ to its corresponding target $\mathbf{y}_n$ within a fixed time $T$. This field is constructed based on purely geometric considerations (see \cref{fig:tubular}) and can be chosen to be Lipschitz continuous.
\begin{proposition}\label{prop:exactsmooth}
Let $N\geq1$, $d\geq2$ and $T>0$ be fixed. Consider the dataset $\mathcal{D}\subset\mathbb{R}^d\times\mathbb{R}^d$ as defined in \eqref{sample}, and any compact subset $\Omega\subset\mathbb{R}^d$ such that $\operatorname{Int}(\Omega)$ is connected and $\mathcal{D}\subset\operatorname{Int}(\Omega)\times\operatorname{Int}(\Omega).$ Then, there exists a vector field $\mathbf{V}\in\operatorname{Lip}\left(\mathbb{R}^d,\mathbb{R}^d\right)$ such that the flow $\Psi_{T,\mathbf{V}}$ of the equation \begin{equation*}\dot\bfx=\mathbf{V}(\bfx)\end{equation*} interpolates the dataset $\mathcal{D}$, and the $N$ curves given by \begin{equation*}\mathcal{C}_n\coloneqq\{\Psi_{t,\mathbf{V}}(\bfx_n):t\in[0,T]\}\qquad (n=1,\dots,N),\end{equation*} are contained in $\operatorname{Int}(\Omega)$.
\end{proposition}

The subset $\Omega\subset\mathbb{R}^d$, which will serve as our domain of approximation, can always be established as $\Omega = [-R, R]^d$ for a sufficiently large $R>0$. Consequently, for any dataset $\mathcal{D} \subset \mathbb{R}^d \times \mathbb{R}^d$, there exists a field $\mathbf{V} \in \operatorname{Lip}(\mathbb{R}^d, \mathbb{R}^d)$ whose integral curves $\mathcal{C}_n$ interpolate $\mathcal{D}$, i.e., the space
\begin{multline*}
\mathcal{V}_\mathcal{D}\coloneqq\big\{\mathbf{V}\in\operatorname{Lip}\left(\mathbb{R}^d,\mathbb{R}^d\right): \Psi_{T,\mathbf{V}} \text{ interpolates }\mathcal{D}
\big\}  
\end{multline*}
is non-empty. Moreover, we can define\begin{equation*}L_0\coloneqq\inf_{\mathbf{V}\in\mathcal{V}_\mathcal{D}}L_V,\end{equation*}
so $L_0$ only depends on $\mathcal{D}$ and the chosen domain of approximation $\Omega$. It suffices then to combine the UAT from \cite{devore_hanin_petrova_2021} (see \cref{lem:devore} in \cref{sec:proofs}) with classical results on the stability of ODEs to obtain the following theorem: 
 \begin{figure}[t!]
   \centering
   \begin{minipage}{0.48\textwidth}
   \centering
    \includegraphics[width=0.98\linewidth]{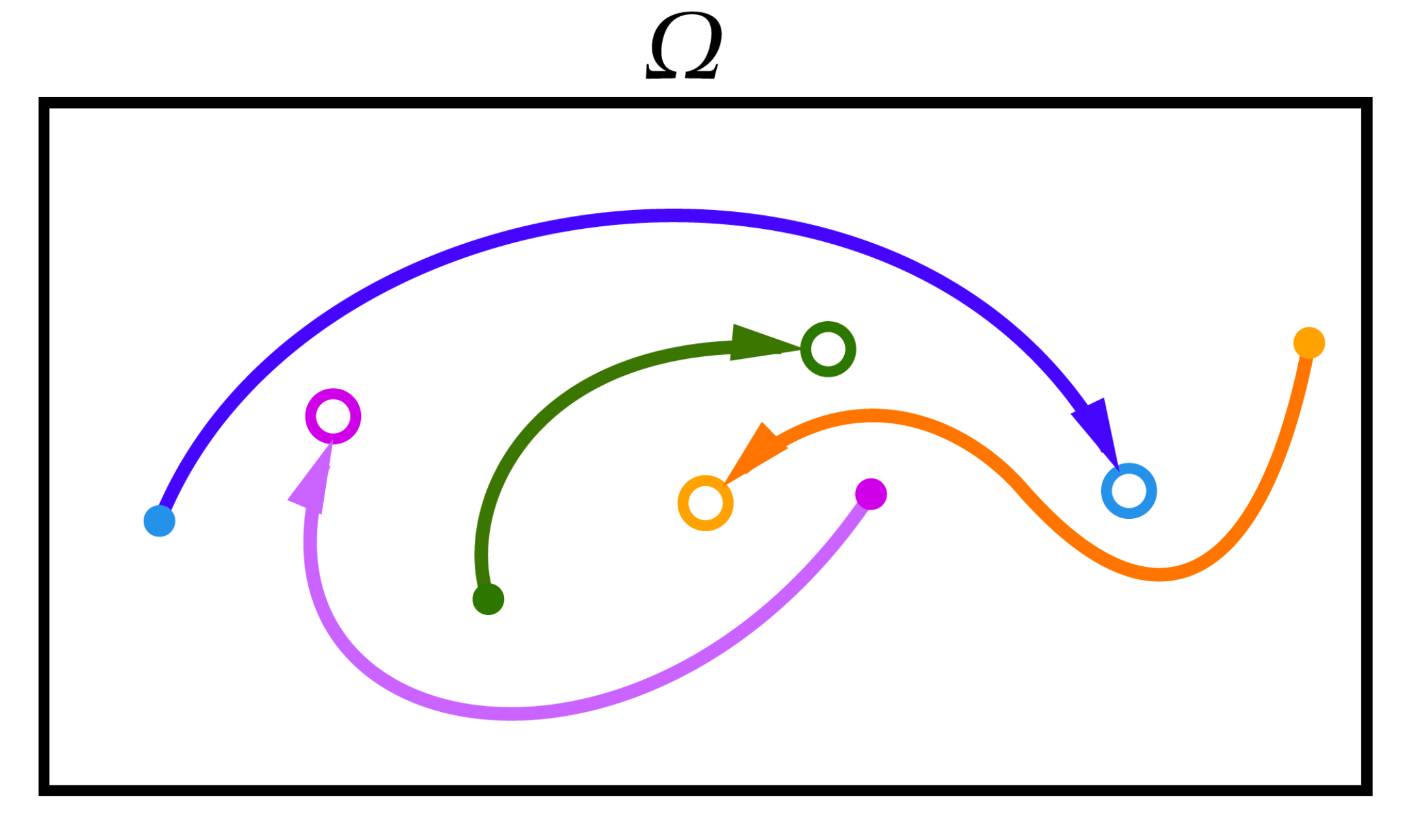}
    \end{minipage}
    \begin{minipage}{0.48\textwidth}
      \centering
    \includegraphics[width=\linewidth]{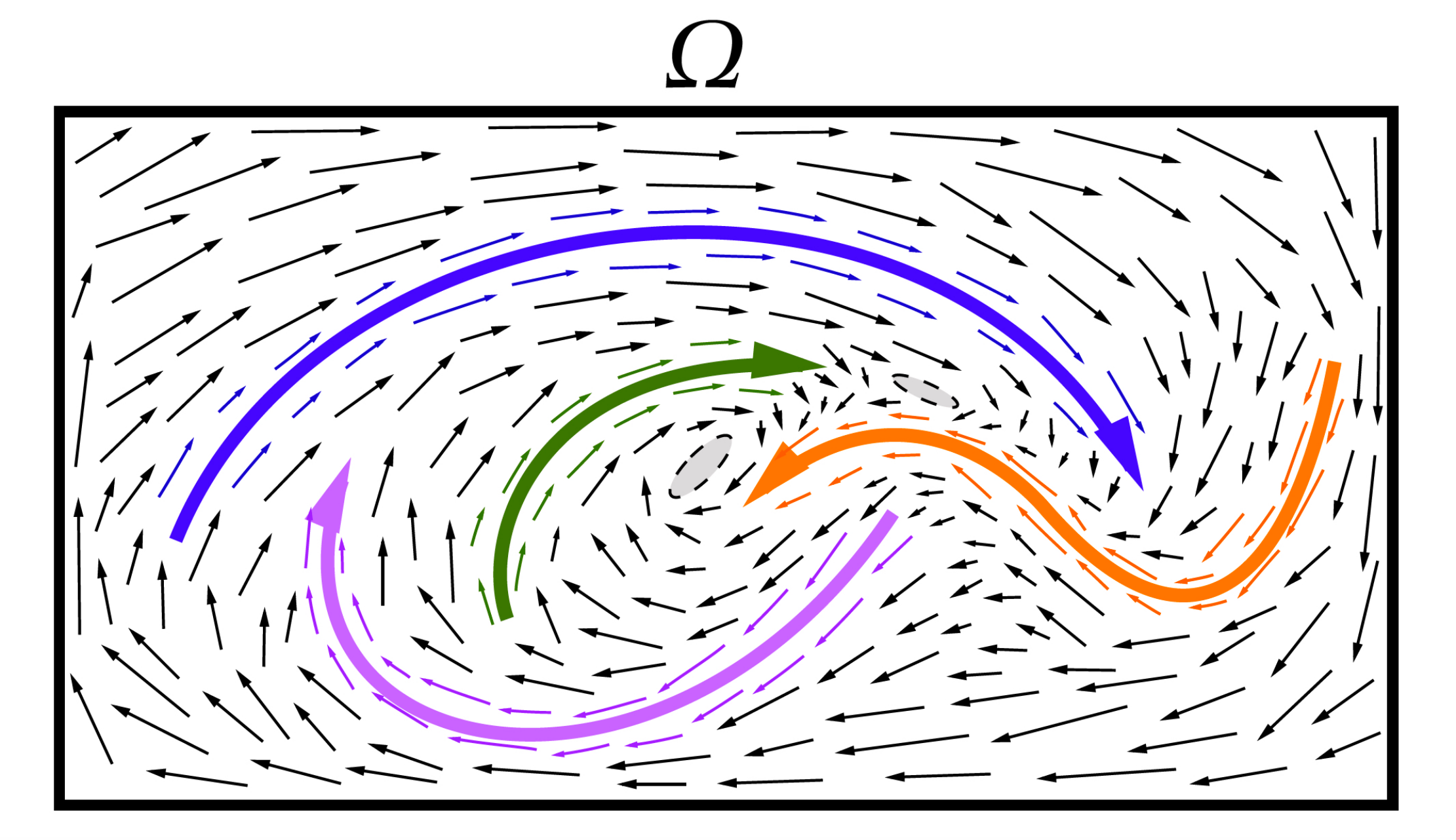}
   \end{minipage}
   \caption{Construction of the Lipschitz field $\mathbf{V}$ in \cref{prop:exactsmooth} which interpolates  $\mathcal{D}$ in a compact domain $\Omega$ that contains all the points and curves.}
\label{fig:tubular}
       \end{figure}
\begin{theorem}\label{th:control-shallow}
 Let $N\geq1$, $d\geq2$ and $T>0$ be fixed. Consider the dataset $\mathcal{D}$ as defined in \eqref{sample}. For each $p\geq1$, there exists a control $(W,A,\mathbf{b})\in \mathbb{R}^{d\times p}\times\mathbb{R}^{p\times d}\times\mathbb{R}^p$ such that the flow $\Phi_T$  generated by \eqref{eq:node-shallow} satisfies
            \begin{equation}\label{eq:ratethm2}
                \sup_{i=1,\dots,N} |\mathbf{y}_n-\Phi_T(\bfx_n;W,A,\mathbf{b})|\leq C_{d,\mathcal{L},T}\,\frac{\log_2(\kappa)}{\kappa^{1/d}},
             \end{equation}
 where $\kappa=(d+2)dp$ is the complexity of the NODE, and
 \begin{equation*}
C_{d,\mathcal{L},T}=C_{d,L_0}T\,\exp\big(\mathcal{L}\, T\big),
 \end{equation*}
 being $\mathcal{L}=\min\{L_0,\|W\|\cdot\|A\|\}$ where $\|\cdot\|$ is the spectral norm, and $C_{d,L_0}>0$ a constant depending on $d$ and $L_0$ but independent of $\kappa$.        
        \end{theorem}
\begin{remark}
The argument employed in this theorem extends beyond neural networks. Since
we solely rely on a density result that provides a convergence rate, other dense families of functions like polynomials, trigonometric, finite element methods or wavelets could also be considered, with their corresponding convergence rates.
\end{remark}

\begin{remark}

Given a domain of approximation $\Omega$, for the bound \eqref{eq:ratethm2} to be optimal, it is natural to pose the problem of finding the interpolating field $\mathbf{V}\in\mathcal{V}_\mathcal{D}$ which has the smallest possible Lipschitz constant $L_V$ within $\Omega$. 
\end{remark}
\begin{remark}
    When $d\geq3$, the construction of a field $\mathbf{V}\in\mathcal{V}_\mathcal{D}$ is generally very simple. Since, in that case, two arbitrary curves are unlikely to intersect, we can generally consider the $N$ segments that connect each pair $(\bfx_n,\mathbf{y}_n)\in\mathcal{D}$ and build $\mathbf{V}$ as one of the piecewise constant fields having these segments as integral curves. Selecting the optimal field then becomes a combinatorial problem.
\end{remark}

\subsection{Transport control}\label{subsec:transp}
   So far, we have considered the system \eqref{eq:node-p} with a finite set of points $\mathcal{D}$ as initial data. A natural extension of this setting, particularly pertinent when the data points are sampled from an underlying distribution, is to consider as input a probability measure $\mu_0$ on $\mathbb{R}^d$. The scenario where this distribution is a finite combination of Dirac deltas aligns with the study previously conducted in \cref{subsec:simcontrol}.

   Specifically, we consider the space $ \mathcal{P}_{ac}^c(\mathbb{R}^d)$ of compactly supported and absolutely continuous probability measures on $\mathbb{R}^d$. Our goal is to transform any given $\mu_0\in\mathcal{P}_{ac}^c(\mathbb{R}^d)$ into a fixed target probability measure $\mu_*$  through the push-forward map generated by a neural ODE, that is, \begin{equation*}\Phi_T(\cdot;W,A,\mathbf{b})_\#\mu_0=\mu_*.\end{equation*}
   This question can be reformulated as the control problem of a transport equation. For each $t\in[0,T]$, we consider the family of measures $\mu(t)=\Phi_{t\#}\mu_0$, where $\Phi_t$ represents the flow at time $t$ generated by \eqref{eq:node-p}. Given that the field \begin{equation*}\sum_{i=1}^p \mathbf{w}_i(t),\sigma(\mathbf{a}_i(t)\cdot\bfx+b_i(t))\end{equation*} is Lipschitz continuous with respect to $\bfx$, if $\mu_0\in\mathcal{P}_{ac}^c(\mathbb{R}^d)$ then the curve of measures $\{\mu(t)\}_{t\in[0,T]}$ is contained in $\mathcal{P}_{ac}^c(\mathbb{R}^d)$. For each $t$, $\mu(t)$ is defined by a density function $\rho(t)$ that satisfies the neural transport equation
    \begin{align}
        \label{eq:NTEq}
        \begin{cases}
             \partial_t \rho +\operatorname{div}_\bfx\left(\rho\sum_{i=1}^p \mathbf{w}_i\,\sigma(\mathbf{a}_i\cdot\bfx+b_i)\right)=0\\
        \rho(0)=\rho_0.
        \end{cases}
    \end{align}
Here, we have assumed that $\mu_0$ has density $\rho_0$,  and \begin{equation*}(\mathbf{w}_i,\mathbf{a}_i,b_i)_{i=1}^p\subset L^\infty\left((0,T);\mathbb{R}^d\times\mathbb{R}^d\times\mathbb{R}\right)\end{equation*}
serve again as control functions.  The projected characteristics of \eqref{eq:NTEq} solve the neural ODE \eqref{eq:node-p}
in $(0,T)\times\mathbb{R}^d$. If the controls are step functions, and since the ReLU function is Lipschitz, the continuity equation \eqref{eq:NTEq} is well-posed and the total mass is conserved. Therefore, we aim to find some controls such that the corresponding solution of \eqref{eq:NTEq} with initial condition $\rho_0$ satisfies    \begin{equation*}\rho(T)=\rho_*.\end{equation*}
    This task, however, can be very hard to achieve, so we consider a relaxation of the problem to approximate control of \eqref{eq:NTEq}. For this purpose, first we must choose a function to quantify the difference between any two measures. 
    \begin{deff}\label{defwasserq}
         For any $q\geq1$, the \emph{Wasserstein-$q$ distance} between $\mu,\nu\in\mathcal{P}_{ac}^c(\mathbb{R}^d)$ is defined as
        \begin{equation}\label{eq:wassdef}       W_q(\mu,\nu)\coloneqq\Big(\min_{\gamma\in \Pi(\mu,\nu)}\int_{\mathbb{R}^d\times\mathbb{R}^d} |\bfx-\mathbf{y}|^qd\gamma(x,y)\Big)^{1/q},
        \end{equation}
        where $\Pi(\mu,\nu)$ denotes the set of measures $\gamma$ on $\mathbb{R}^d\times\mathbb{R}^d$ that couple $\mu$ and $\nu$ in the sense that $\gamma(\cdot\times\mathbb{R}^d)=\mu(\cdot)$ and $\gamma(\mathbb{R}^d\times \cdot)=\nu(\cdot)$. Note that $\mu\in\mathcal{P}_{ac}^c(\mathbb{R}^d)$ has finite $q$-th momentum for every $q\geq1$, hence the Wasserstein$-q$ distance is well-defined in this space. Moreover, recalling the Monge formulation of optimal transport, if $\mu$ and $\nu$ belong to $\mathcal{P}_{ac}^c(\mathbb{R}^d)$ then 
                \begin{equation}\label{eq:monge} W_q(\mu,\nu)=\Big(\min_{T}\Big\{\int_{\mathbb{R}^d}|\bfx-T(\bfx)|^q d\mu: T_\#\mu=\nu\Big\}\Big)^{1/q},
                \end{equation}
                where $T:\mathbb{R}^d\rightarrow\mathbb{R}^d$ measurable, see \cite{VILLANI_2016}.
    \end{deff}

     \begin{probdef*}
     Let $\mu_0$ and $\mu_*$ be two compactly supported, absolutely continuous probability measures with respective densities $\rho_0$ and $\rho_*$. For any fixed time horizon $T>0$ and $\varepsilon>0$, find controls \begin{equation*}\left\{(\mathbf{w}_i,\mathbf{a}_i,b_i)\right\}_{i=1}^p\subset L^\infty\left((0,T);\mathbb{R}^d\times\mathbb{R}^d\times\mathbb{R}\right),\end{equation*}for some $p\geq1$, such that the solution of \eqref{eq:NTEq} in time $T$ \emph{approximately interpolates} the initial condition $\rho_0$ to the target density $\rho_*$. This is achieved when the $W_q$-error of the corresponding measures (for some $q\geq1$) satisfies:
    \begin{equation*}W_q(\mu(T),\mu_*)<\varepsilon.\end{equation*}
 \end{probdef*}  
The following theorem offers a partial solution to this problem. It assumes that $1\leq q< \frac{d}{d-1}$ and targets the uniform measure in $[0,1]^d$. While this bears resemblance to achieving null controllability, the nonlinear nature of the problem prevents from directly extending this result to arbitrary targets.

\begin{theorem}\label{th:NTEq-control-p}
         Let $d\geq1$, $\mu_0\in\mathcal{P}_{ac}^c(\mathbb{R}^d)$ with density $\rho_0$, $\mu_*$ the uniform measure in $[0,1]^d$, and $T>0$ be fixed.  For any $\varepsilon>0$, $q\in[1,\frac{d}{d-1})$ and $p\geq1$, there exists a piecewise constant control 
         \begin{equation*}\left(W,A,\mathbf{b}\right)\in L^\infty\left((0,T);\mathbb{R}^{d\times p}\times\mathbb{R}^{p\times d}\times\mathbb{R}^p\right)\end{equation*}
    such that the measure $\mu(t)\in\mathcal{P}_{ac}^c(\mathbb{R}^d)$ whose density $\rho(t)$ solves \eqref{eq:neuraltransport} taking $\rho_0$ as initial condition, satisfies
 \begin{equation*}W_q(\mu(T),\mu_*) < \varepsilon.\end{equation*}
    Furthermore,  the number of discontinuities of $(W,A,\mathbf{b})$ is \begin{equation*}L=\left\lceil 2d/p\right\rceil+\max\{ \lceil n/p_1\rceil,\dots, \lceil n^d/p_d\rceil\}-1, \end{equation*}
    for any $p_1,\dots,p_d\geq1$ such that $p_1+\cdots+p_d=p$, and \begin{equation*}       n\coloneqq\left(\frac{3 d^{1/2+1/q}}{\varepsilon}\right)^{\frac{1}{1+d/q-d}}.
    \end{equation*}
        \end{theorem}  
For a given $\varepsilon>0$, the behavior of $L$ resembles that described in \cref{th:SimControl-p}, as it decreases with an increase in $p$, reaching $L=1$ when $p$ is large enough. Our proof is similar to a strategy from \cite{velocity-control}, and based on the specific movements that the neural ODE \eqref{eq:node-p} allows. We compress the support of $\mu_0$ to $[0,1]^d$ and divide it into hyperrectangles, each with a mass of $O(\varepsilon^d)$. These subsets are then transformed to match a similar partition of $[0,1]^d$ corresponding to the uniform measure $\mu_*$.
       
    \begin{remark}
        If $\varepsilon>0$ is sufficiently small, and we choose $p_1=\cdots=p_{d-1}=1$, $p_d=p-d+1$, it follows that
        \begin{equation*}
            L=\lceil 2d/p\rceil +\left\lceil\frac{1}{p-d+1} \left(\frac{3^{1+d/q}\sqrt{d}}{\varepsilon}\right)^{\frac{d}{1+d/q-d}}\right\rceil-1.
        \end{equation*}
      For $q=1$, this expression simplifies to:
           \begin{equation*}
               L=\left\lceil 2d/p\right\rceil+\left\lceil \frac{1}{p-d+1}\left(\frac{3^{1+d}\sqrt{d}}{\varepsilon}\right)^d\right\rceil-1.
           \end{equation*}
    \end{remark}

\section{Discussion}\label{conclusions}
\subsection{Conclusions}
We have established several results on the capacity of neural ODEs for interpolation and its relationship with the chosen architecture, determined by the depth $p$ and width $L$. More precisely, we have provided explicit dependencies between these two parameters that are sufficient to (exactly or approximately) interpolate either two sets of $N$ different points in $\mathbb{R}^d$ or any compactly supported, absolutely continuous probability measure in $\mathbb{R}^d$ with the uniform measure in $[0,1]^d$. Our work reveals that $p$ and $L$ can play similar roles in the algorithms, thereby exhibiting a degree of exchangeability in the network's structure.

Specifically, \cref{th:SimControl-p} proves that a neural ODE with $p$ neurons can interpolate any dataset of $N$ pairs of points using a piecewise constant control with $L = 2\left\lceil N/p\right\rceil - 1$ discontinuities.  Although increasing $p$ reduces the number of discontinuities, we find a limiting case of a 2-hidden layer neural ODE when $p\geq N$. Explicit controls for interpolation with a shallow neural ODE $(L=0)$ are obtained in \cref{cor:change-basis} when $d>N$; or in \cref{cor:snodewithsep}, with $p=N$ under \cref{hyp1}. More generally, \cref{th:control-shallow} provides an error decay rate with respect to the number of parameters for shallow neural ODEs.  Finally, \cref{th:NTEq-control-p} explores the Wasserstein-$q$ approximate control of the neural transport equation to a uniform distribution on $[0,1]^d$ using piecewise constant controls. As in \cref{th:SimControl-p}, the number of discontinuities diminishes as $p$ increases. 
\subsection{Open questions}
Some new objectives can be derived from our work:\newline
\textbf{1. Approaching the autonomous regime.} 
        As we have discussed, both $d>N$ and \cref{hyp1} are only special cases where we can find controls to interpolate in the autonomous regime of shallow neural ODE. The question of finding such a construction for any $d\geq1$, or at least under a less restrictive hypothesis on the dataset than \cref{hyp1}, is still open. A first step could involve assuming the relaxed condition that the projections of points onto a line with direction $\mathbf{a}\in\mathbb{S}^{d-1}$ are ordered as $\mathbf{a}\cdot\bfx_{\tau(1)} < \cdots < \mathbf{a} \cdot \bfx_{\tau(N)}$ and $\mathbf{a}\cdot \mathbf{y}_{\tau(1)} < \cdots < \mathbf{a} \cdot \mathbf{y}_{\tau(N)}$, for a certain permutation $\tau$ of $N$ elements. The strategy to be adopted is clear. It entails a combination of the one-dimensional control delineated in \cref{lemma:control-subdiv} with the transversal control required to prove \cref{cor:snodewithsep}, followed by a case-by-case analysis.
\newline
\textbf{2. Universal approximation.}
  In \cite{rbz-node}, the authors demonstrate that a narrow neural ODE can approximate any simple function with compact support. The proof hinges on three key aspects: interpolation capacity, the compressive nature of neural ODEs, and a control strategy that bounds the support of the target function irrespective of required approximation accuracy. However, minimizing the number of time discontinuities required for this control, particularly by increasing $p$, is a non-trivial task that may require an entirely different approach.\newline
  \textbf{3. Neural transport equation.}
        The problem of controlling the neural transport equation \eqref{eq:NTEq} using constant controls is as yet unresolved. One possible approach involves approximating both the initial and target measures with atomic measures of the form $\rho_N = \frac{1}{N} \sum_{n=1}^N \beta_n\delta_{\alpha_n}$, with $\beta_n>0$ and $\alpha_n\in\mathbb{R}^d$, and then interpolating those Dirac deltas by controlling the characteristic curves. However, a potential issue arises as $N \to \infty$: the distance between $\rho(T)$ and $\rho^N(T)$, the solutions to the transport equation with initial conditions $\rho_0$ and $\rho_0^N$, respectively, may diverge significantly. This error can be quantified using the Grönwall inequality, which suggests that the Lipschitz constant could increase unboundedly with $N$, especially as the number of controlled points grows.\newline
        \textbf{4. Minimizing the number of time jumps.}
        Another interesting question is how to frame the reduction of discontinuities as an optimal control problem. For instance, one could penalize the frequency of time jumps by targeting the total variation seminorm. However, this seminorm lacks regularity, and moreover the class of piecewise constant functions is not a closed set of admissible controls.\newline
 \textbf{5. Switching dimensions.}  
In our simplified ResNet \eqref{eq:resnet-intro}, the dimension remains constant across layers. However, strategically varying the hidden dimension by defining $p=p(t)$ could offer advantages, either by reducing complexity through dimension shrinkage or by creating space through dimension increase. Exploring effective methods to implement these transitions, whether by employing projections or by applying nonlinear transformations to the data, constitutes a compelling area for research.
\section{Proofs}\label{sec:proofs}
\subsection{Basic dynamics}
\begin{figure}
 \includegraphics[scale=0.049]{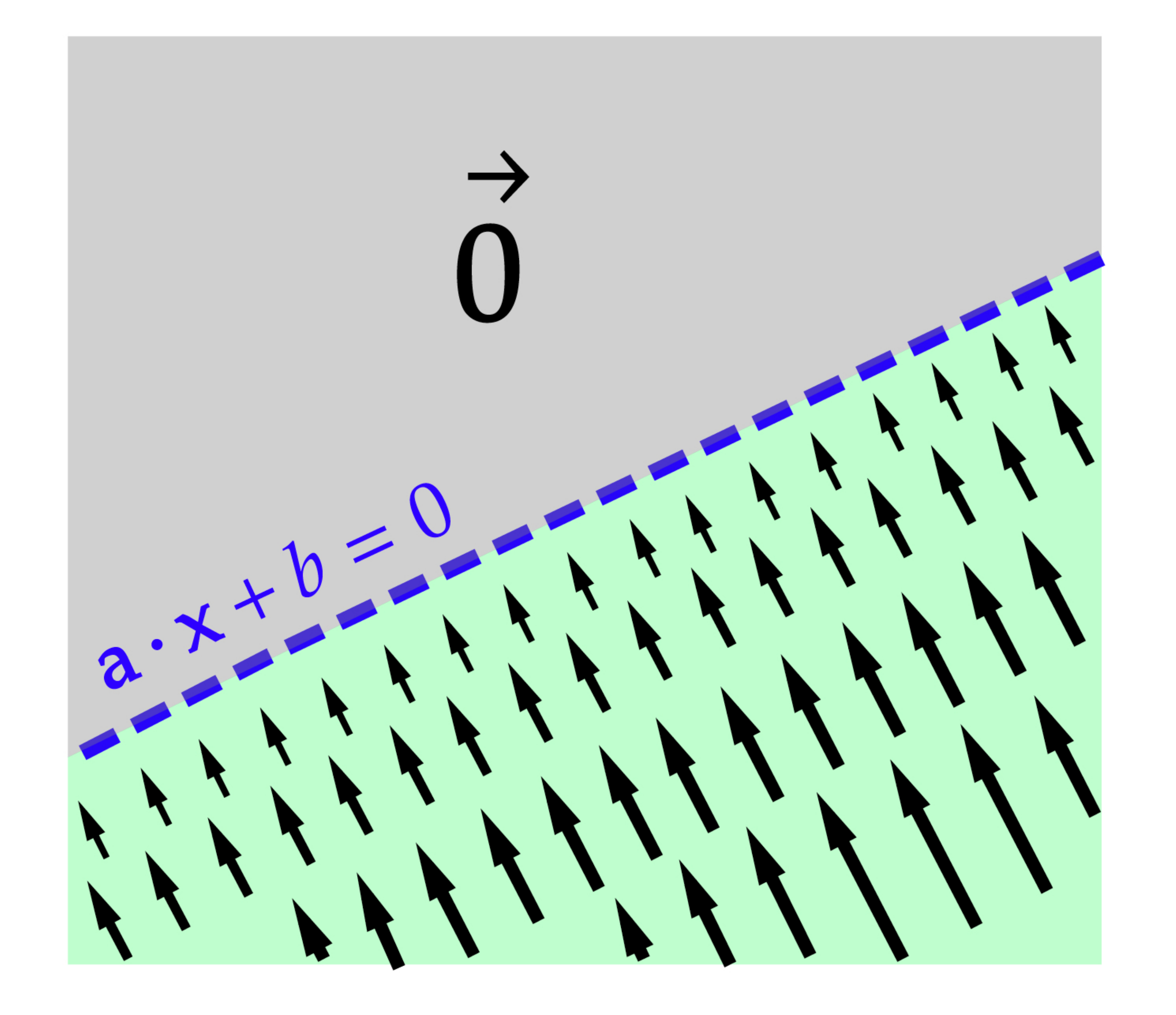}
  \includegraphics[scale=0.05]{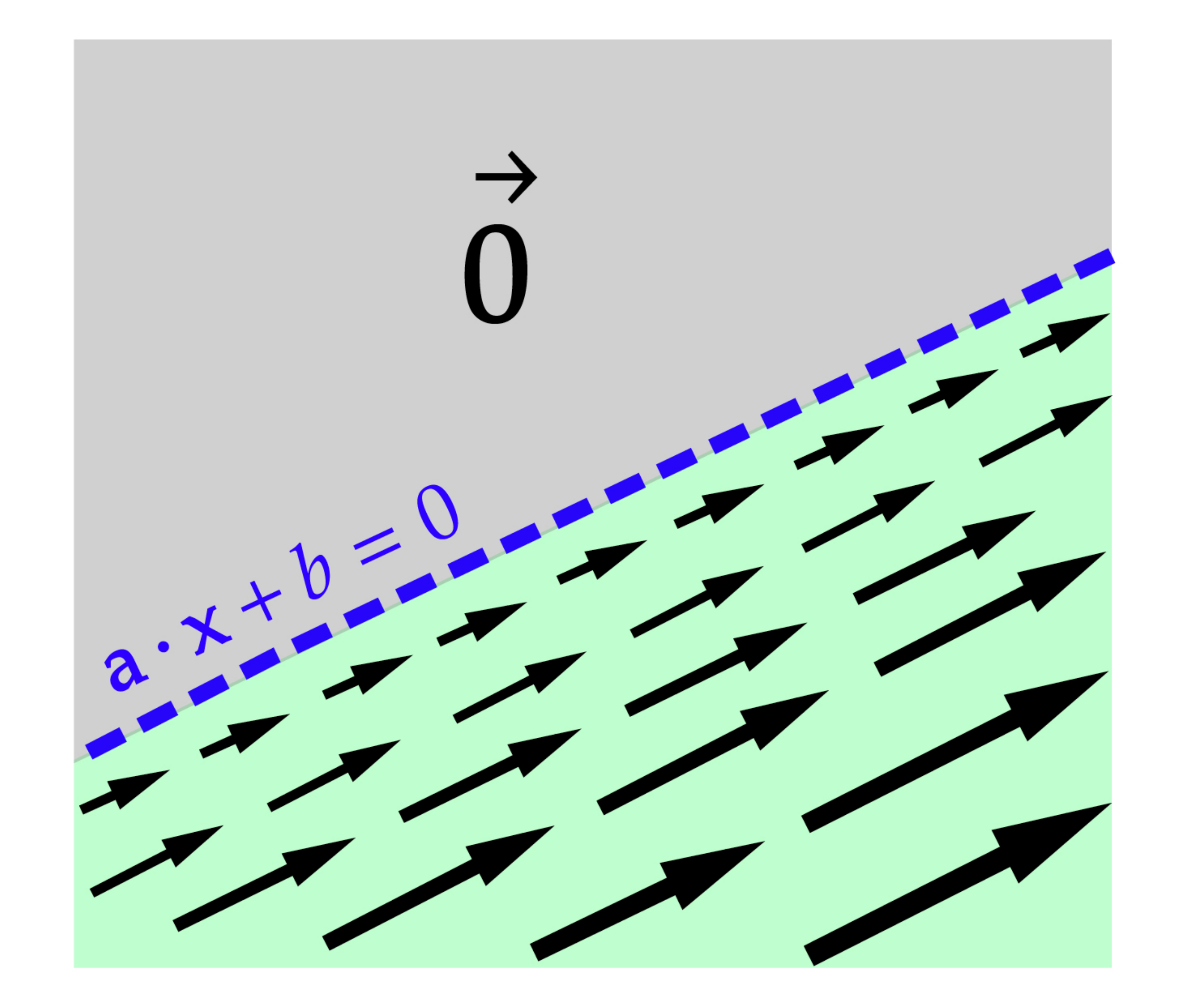}
 \includegraphics[scale=0.05]{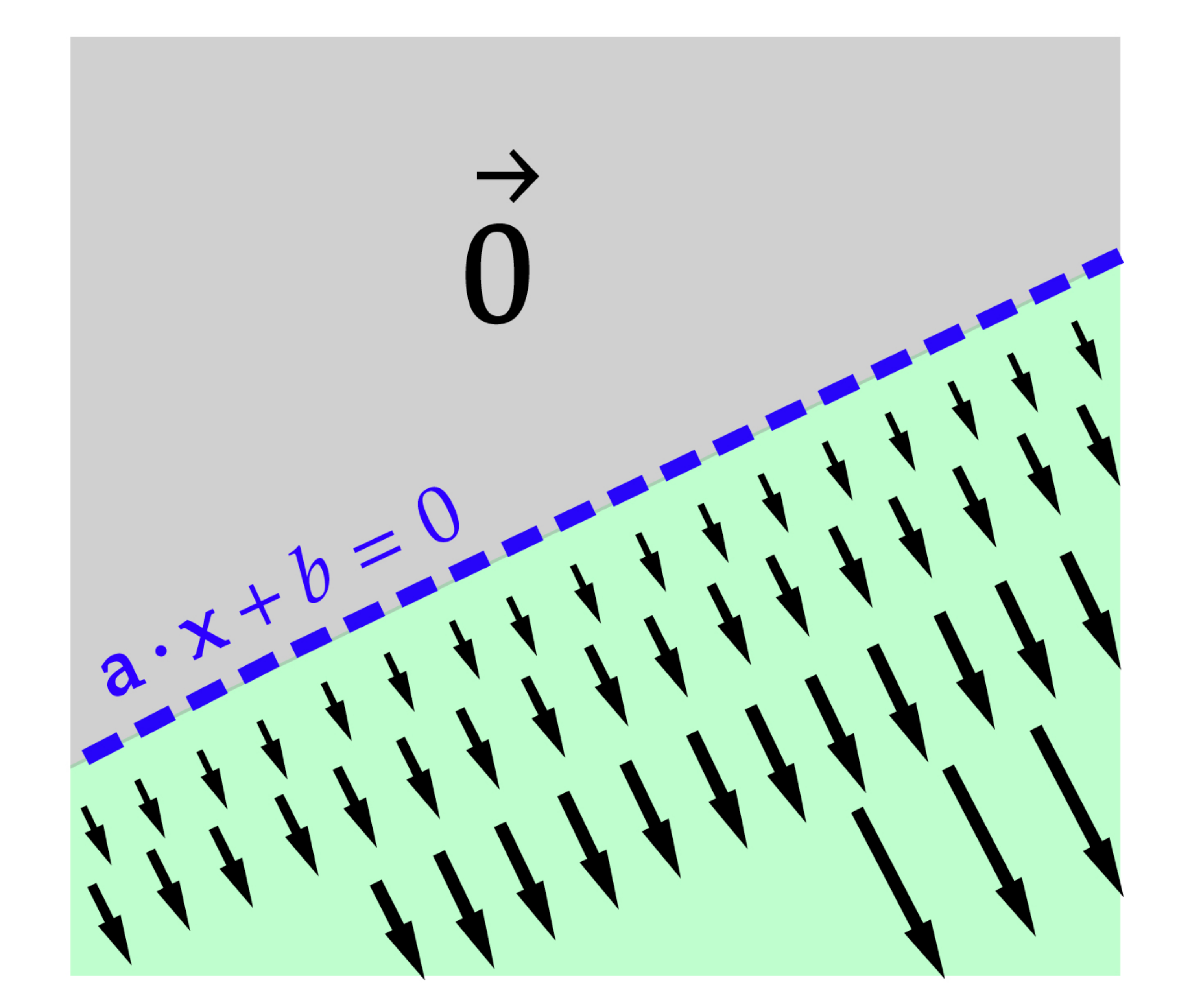}
 \caption{Left to right: Compression, parallel motion, expansion.}
 \label{fig:basicmove}
\end{figure}
We describe the simplest dynamics that we can generate via \eqref{eq:node-Narrow} by conveniently choosing $(\mathbf{w},\mathbf{a},b)$:

1. For each $t>0$, the term $\mathbf{a}(t)\cdot\bfx+b(t)$ identifies a hyperplane in $\mathbb{R}^d$. For instance, taking $\mathbf{a}= \mathbf{e}_k$ and $b= -c$, we fix the hyperplane $h$ with the equation $x^{(k)}-c=0$.

2. The application of $\sigma$ and the product with the vector $\mathbf{w}(t)$ yields the field $\mathbf{w}(t)\max\{x^{(k)}-c,0\}$, which exhibits distinct dynamics in two complementary half-spaces: $H^+ \equiv \{x^{(k)}-c>0\}$, where the field equals $\mathbf{w}(t)(x^{(k)}-c)$, and $H^- \equiv \{x^{(k)}-c\leq0\}$, where the field is zero, meaning this set remains stationary under the flow.

3. The choice of $\mathbf{w}(t)$ specifies the orientation and magnitude of the field. For example, $\mathbf{w}(t)=\pm\mathbf{e}_k$ yields $\dot{x}^{(k)}(t)= \pm\max\{(x^{(k)}-c),0\}$, so the points in $H^+$ can either be attracted to or repelled from $h$, enabling compression or expansion along the $k$-th coordinate. Conversely, $\mathbf{w}(t)=\mathbf{e}_i$, with $i\not=k$, results in $\dot{x}^{(i)}(t)=\max\{(x^{(k)}-c),0\}$. In this case, we generate in $H^+$ a movement that is parallel to $h$, i.e., along the coordinate $i$.

    The three basic operations of compression, expansion and movement in parallel with the hyperplane (represented in \cref{fig:basicmove}) constitute our toolbox for many subsequent proofs.        
   
    \subsection{Proof of \cref{th:SimControl-p}.}\label{appendix:simControl-p}
\begin{figure}
\centering
\begin{subfigure}{0.45\linewidth}
    \includegraphics[width=0.95\linewidth]{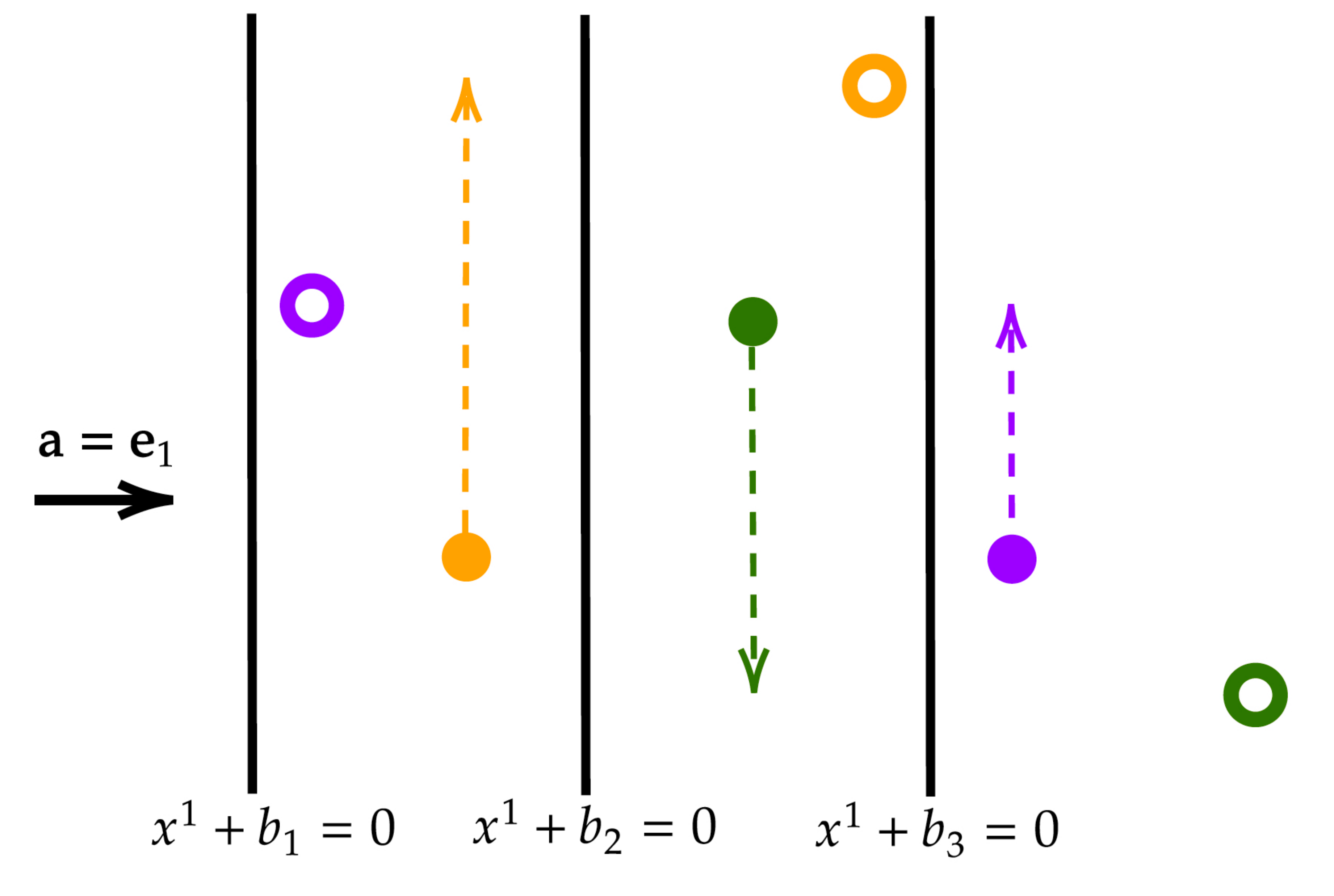}
    \label{fig:node-control1}
\end{subfigure}
\begin{subfigure}{0.45\linewidth}
    \includegraphics[width=0.95\linewidth]{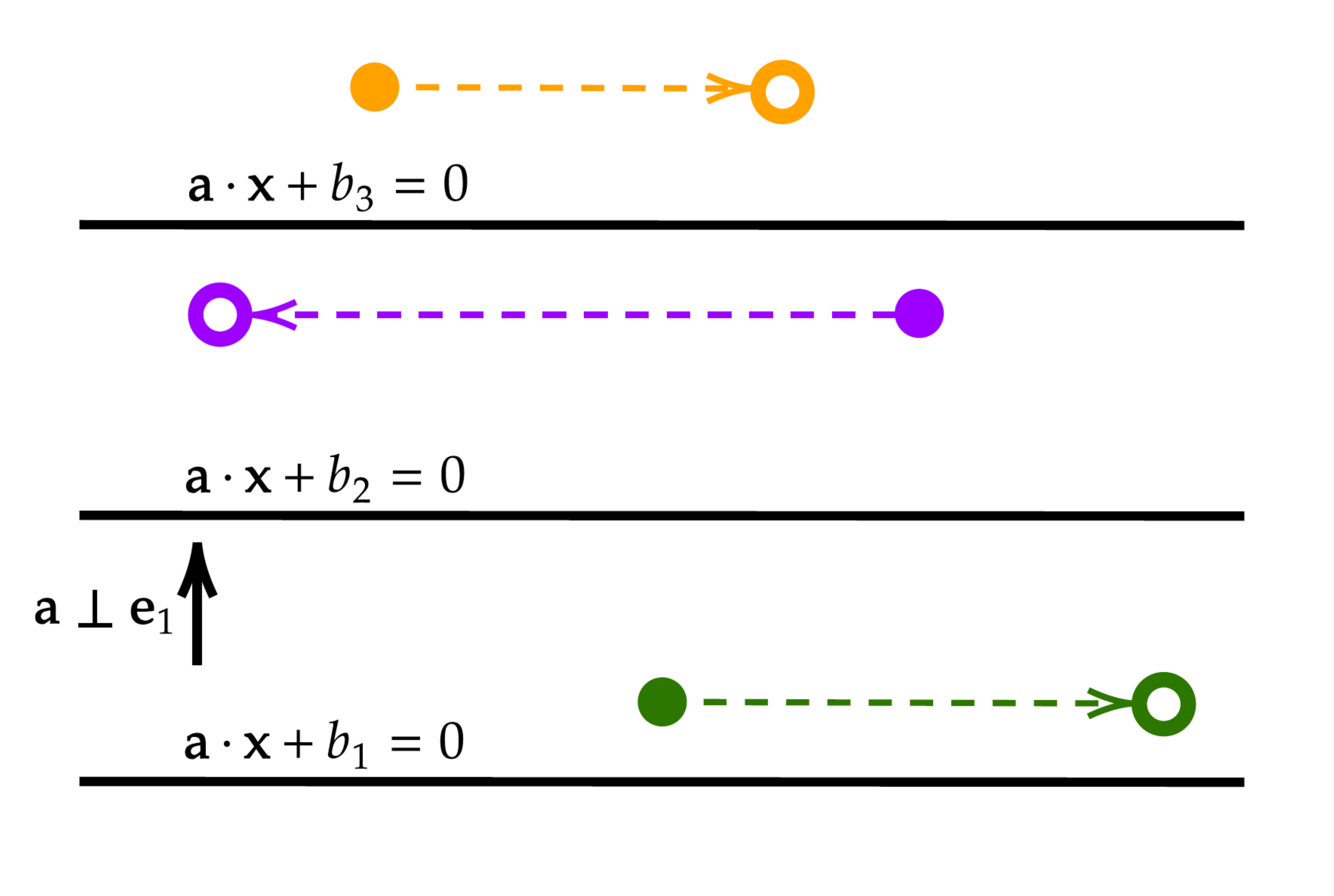}
    \label{fig:node-control2}
\end{subfigure}
\caption{Left: Step 1. Fix $x^{(1)}$ and control $x^{(2)},\dots,x^{(d)}$. Right: Step 2. Control $x^{(1)}$ while $x^{(2)},\dots,x^{(d)}$ are fixed.}
\label{fig:node-control}
\end{figure}

    We will employ the following lemma, whose proof we postpone to the end of this subsection:
    \begin{lemma}\label{lemma:SimControl-p} Let $N\geq1$, $d\geq2$ and consider the dataset $\mathcal{D}=\{(\bfx_n,\mathbf{y}_n)\}_{n=1}^N$ as defined in \eqref{sample}.  There exists a change of coordinates in $\mathbb{R}^d$ such that \begin{equation}\label{hyp:separation-SimControl}x_n^{(1)}\not=x_m^{(1)}\qquad\text{and}\qquad y_n^{(2)}\not=y_m^{(2)}, \qquad \text{if } n\not=m.             \end{equation} 
    \end{lemma}
  Under the separability condition \eqref{hyp:separation-SimControl}, we achieve the exact control by building on the methods developed in \cite{rbz-node}. Let $p\geq1$ be fixed. We divide the proof in two steps, illustrated in \cref{fig:node-control1} and \cref{fig:node-control2}.\newline
 \textbf{Step 1: Control of $d-1$ coordinates.}
By \eqref{hyp:separation-SimControl}, we can relabel the data $\{\mathbf{x}_n\}_{n=1}^N$ to impose the ordering \begin{equation*}x_1^{(1)} < \cdots < x_{N}^{(1)}.\end{equation*} We define a partition of $\{\mathbf{x}_n\}_{n=1}^N$ in $\lceil N/p\rceil$ subsets by increasing order of the $x^{(1)}$-coordinate. The $j$-th subset is
\begin{equation*}X_j\coloneqq\{\bfx_{(j-1)\cdot p+1},\dots,\bfx_{j\cdot p}\},\quad\text{for } j=1,\dots,\lceil N/p\rceil-1,\end{equation*} and $X_{\lceil N/p\rceil}$ contains the remaining $N-p\lfloor N/p\rfloor$ points. 
We describe the control of the first subset $X_1$. We take controls $\mathbf{a}_i=\mathbf{e}_1$ and $b_i\in\mathbb{R}$, for $i=1,\dots,p$, satisfying \begin{equation*}-b_1<x_{1}^{(1)}<-b_2<x_{2}^{(1)}<\cdots <-b_p<x_{p}^{(1)}.\end{equation*}
These controls define a family of parallel hyperplanes, given by    $\mathbf{a}_i\cdot\bfx+b_i=x^{(1)}+b_i=0$,
which separate the points of $X$. In this way, $i-1$ terms of the sum in \eqref{eq:node-p} cancel inside the strip $-b_i<x^{(1)}<-b_{i+1}$ for each $i=1,\dots,p$, so \eqref{eq:node-p} simplifies to
\begin{equation}\label{nodestripn}
\dot\bfx=\sum_{l=1}^i\big\{\mathbf{w}_lx^{(1)}+\mathbf{w}_lb_l\big\}.
\end{equation} 
We consider velocities of the form $\mathbf{w}_i=(0,w_i^{(2)},\dots,w_i^{(d)}),$ where the components $w_i^{(k)}\in\mathbb{R}$ have to be defined in order to achieve the exact control in time $T=1$.  The first point, $\bfx_{1}$, is subject only to one velocity, $\mathbf{w}_1$, so \begin{equation*}x_{1}^{(k)}(t)=w_{1}^{(k)}(x_{1}^{(1)}+b_1)t+x_{1}^{(k)},\end{equation*}
               while $x_{1}^{(1)}$ remains fixed.  Therefore, it is enough to take 
             \begin{equation*}w_{1}^{(k)}= \frac{y_{1}^{(k)}-x_{1}^{(k)}}{x_{1}^{(1)}+b_1}.\end{equation*}
            Similarly, for $i=2,\dots,p$, having fixed $\mathbf{w}_1,\dots,\mathbf{w}_{i-1}$ it is enough to take 
            \begin{equation*}w_i^{(k)}=\frac{y_{i}^{(k)}-x_{i}^{(k)}-\sum_{l=1}^{i-1}w_{l}^{(k)}(x_{i}^{(1)}+b_l)}{x_{i}^{(1)}+b_i},\end{equation*}
            for $k=1,\dots,d$. The described procedure can be simultaneously done for each $X_j$, with $j=2,\dots\lceil N/p\rceil$, taking into account that the fields used to control $X_1,\dots,X_{j-1}$ (all of them orthogonal to $\mathbf{e}_1$) will be added as new terms in \eqref{nodestripn}. In the end, we will have, for every $n=1,\dots,N$:
            \begin{equation*}
\Phi_1(\bfx_{n};W,A,\mathbf{b})^{(k)}=y_{n}^{(k)},\qquad\text{for } k=2,\dots,d.
            \end{equation*}
        The total number of iterations employed in this step is $\left\lceil N/p \right\rceil$, which corresponds to $\left\lceil N/p\right\rceil-1$ switches.\newline
         \textbf{Step 2: Control of the remaining coordinate.}  
           In a slight abuse of notation, we redefine $\bfx_n \coloneqq \Phi_1(\bfx_n)$, where $\Phi_1$ is the flow resulting from step 1. Once again, we can relabel the data, now assuming \begin{equation*}x_1^{(2)} < \dots < x_N^{(2)}.\end{equation*} Following the increasing order of the $x^{(2)}$-coordinate, we define $X_1, \dots, X_{\lceil N/p \rceil-1}$, each being a subset of $\{\mathbf{x}_n\}_{n=1}^N$ with $p$ points, and $X_{\lceil N/p \rceil}$, which contains the remaining $N-p\lfloor N/p \rfloor$ points.

We follow an analogous methodology to step 1. For each $j$, we define controls $\mathbf{a}_i = \mathbf{e}_2$ and $b_i$ (for $i=1, \dots, p$) that separate the points of $X_j$ using $p$ parallel hyperplanes, each described by the equation $x^{(2)} = b_i$. Now, we consider velocities of the form $\mathbf{w}_i = w_i\mathbf{e}_1$, where the values $w_i$ are determined, as in step 1, to ensure
            \begin{equation*}
\Phi_1(\bfx_{n};W,A,\mathbf{b})^{(1)}=y_{n}^{(1)},\qquad\text{for } n=1,\dots,N.
            \end{equation*}
          The number of switches employed in step 2 is $\left\lceil N/p\right\rceil-1.$ Then, by adding one more to transition between steps, the whole control requires $L=2\left\lceil N/p\right\rceil-1$ switches, hence proving \cref{th:SimControl-p}.
      \begin{proof}[Proof of \cref{lemma:SimControl-p}]\label{proof:simControl}
      The set of vectors in $\mathbb{R}^d$ that are orthogonal to any point $\bfx_n$ from $\mathcal{D}$, is a finite union of hyperplanes in $\mathbb{R}^d$. Therefore, one can always choose $\mathbf{u}_1\in\mathbb{S}^{d-1}$ inside the complement of this set in $\mathbb{R}^d$. 
      
      Consider the orthogonal subspace $S=\langle\mathbf{u}_1\rangle^\perp\subset\mathbb{R}^d$, of dimension $d-1$. With a similar argument, one can select a vector $\mathbf{u}_2\in\mathbb{S}^{d-1}\cap S$ that has a non-zero scalar product with all of the vectors $\mathbf{y}_n$ $(n=1,\dots,N)$.      
       Completing the pair $(\mathbf{u}_1,\mathbf{u}_2)$ to form a orthonormal basis $\mathcal{B}=\{\mathbf{u}_1,\dots,\mathbf{u}_d\}$ of $\mathbb{R}^d$, the dataset $\mathcal{D}$ will satisfy the separability condition \cref{hyp:separation-SimControl} when expressed in $\mathcal{B}$.
          \end{proof}
\subsection{Proof of \cref{cor:change-basis}}
We aim to eliminate the initial step in the algorithm defined in the proof of \cref{th:SimControl-p}. To achieve this, we seek a new vector basis in $\mathbb{R}^d$ where the input-target pairs inherently share the first coordinate. When $d\leq N$, this condition usually cannot be met. However, when $d>N$, there exists an orthonormal vector basis $\mathcal{B}\subset\mathbb{R}^d$ such that for all $n = 1, \ldots, N$, the first coordinates of the $N$ pairs with respect to $\mathcal{B}$ satisfy $x_n^{(1)} = y_n^{(1)}$.

To construct such a vector system, without loss of generality we can assume that $d=N+1$. Let $(\bfx,\mathbf{y})\in\mathbb{R}^d\times\mathbb{R}^d$ with $\bfx\neq\mathbf{y}$. We seek a vector $\mathbf{u}\in\mathbb{S}^{d-1}$ such that $\mathbf{u}\cdot\bfx=\mathbf{u}\cdot\mathbf{y}$. This condition is equivalent to $\mathbf{u}\cdot(\bfx-\mathbf{y})=0$, which is satisfied by any unit vector $\mathbf{u}$ contained in the linear hyperplane orthogonal to $\bfx-\mathbf{y}$. 

\begin{figure}[t]
   \centering
   \begin{minipage}{0.49\textwidth}
   \centering
    \includegraphics[width=0.62\linewidth]{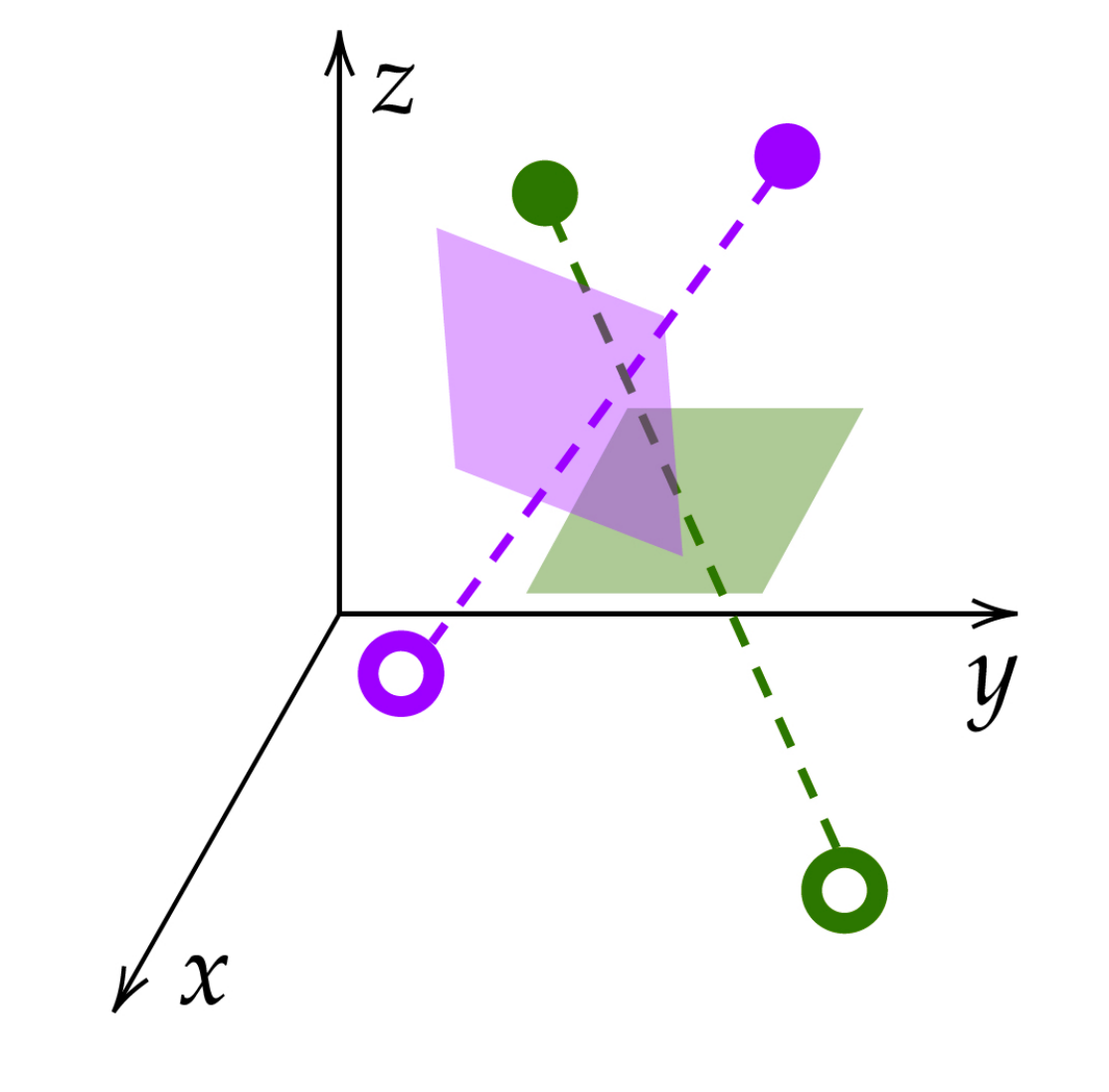}
    \end{minipage}
    \begin{minipage}{0.49\textwidth}
      \centering
    \includegraphics[width=0.9\linewidth]{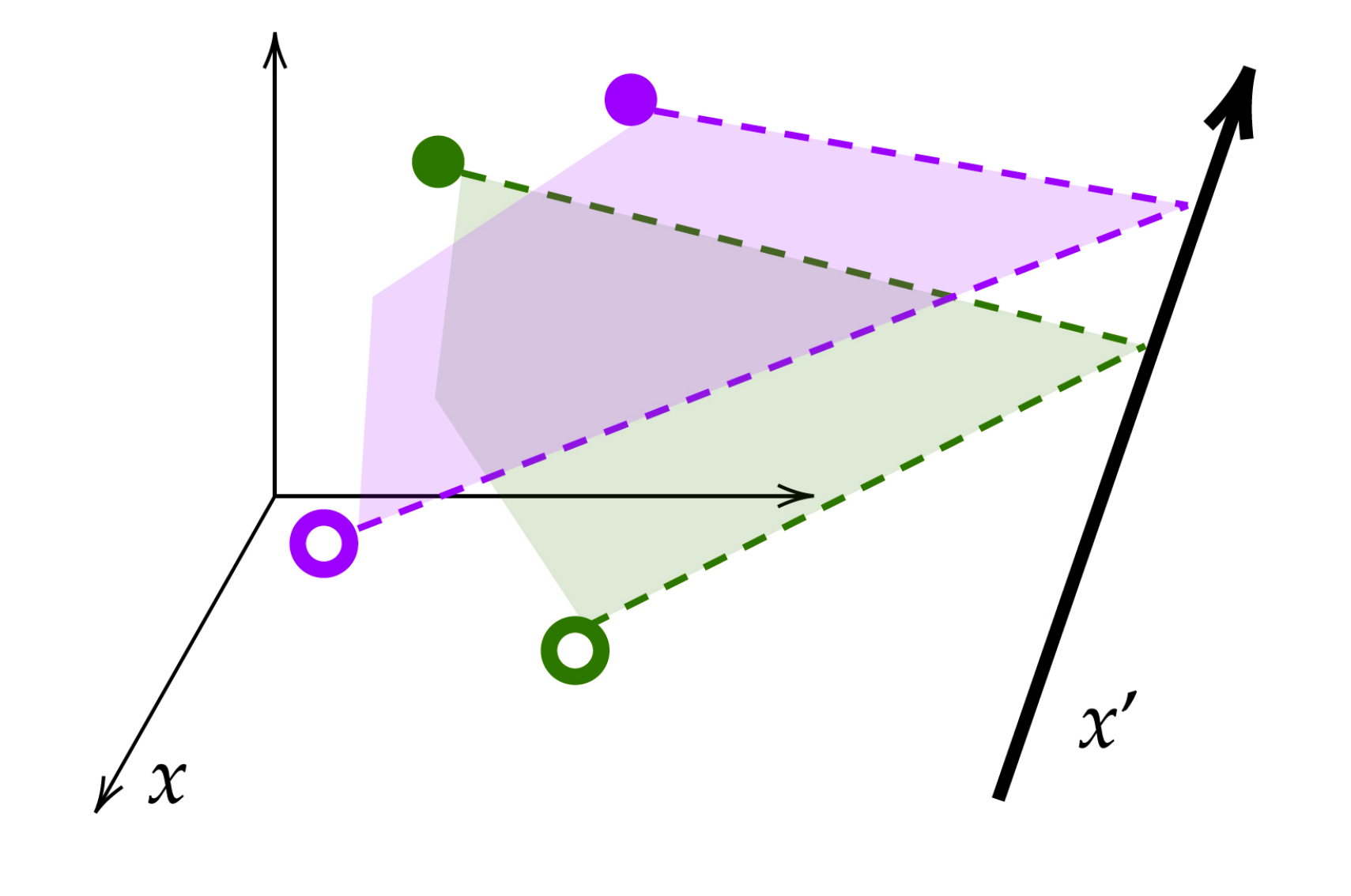}
   \end{minipage}
   \caption{For $N=2$, $d=3$, construction of a new basis of $\mathbb{R}^d$ in which the first coordinates of the pairs $(\bfx_n,\mathbf{y}_n)$ are matched. \cref{th:SimControl-p} is applied afterwards.}
\label{fig:dgtrn}
       \end{figure}

For $d-1$ input-target pairs of points, we consider the corresponding hyperplanes $\{H_{n}\}_{n=1}^{d-1}$. Note that some of these hyperplanes can be repeated. So, the intersection $\bigcap_{i=1}^{d-1} H_i$ yields a linear subspace of dimension at most $d-1$. We choose any unit vector $\mathbf{e}_1'$ contained in that subspace. Any completion to an orthonormal basis $\mathcal{B}=\{\mathbf{e}_1',\dots,\mathbf{e}_d\}\subset\mathbb{R}^d$ will satisfy the desired condition for $d-1$ points. The procedure is illustrated in \cref{fig:dgtrn}.
           
           \subsection{Proof of \cref{cor:snodewithsep}}
     \begin{proof}[Proof of \cref{lem:probability}]
We carry out a similar analysis to the one in \cite{alvlop24}, where the probability of requiring $k$ hyperplanes to separate two sets of $N$ points for binary classification was estimated. We introduce the random variable $Z^p_{d,2N}(\mathcal{D})$, which assigns to each possible dataset \begin{equation*}\mathcal{D}=\{(\bfx_n,\mathbf{y}_n)\}_{n=1}^N\subset\operatorname{supp}(\mu)\times\operatorname{supp}(\mu)\end{equation*}the minimum number of parallel hyperplanes needed to separate in $\mathbb{R}^d$ every pair $(\bfx_n,\mathbf{y}_n)$ from the others. Note that $\min_{\mathcal{D}}Z^p_{d,2N}(\mathcal{D})= N-1$. We estimate the probability $P_{d,N}\coloneqq P(Z_{d,2N}^p(\mathcal{D})=N-1)$ for any $d$ and $N$. 

First, we consider the one-dimensional case. Since all the points are sampled from the same distribution, every possible configuration of the $2N$ points in the real line will have the same probability, i.e., their distribution is uniform on the finite space of all possible orderings. Therefore, we can compute: 
\begin{equation*}
P_{1,N}=\frac{\text{favorable configurations}}{\text{total configurations}}=\frac{N!\,2^N}{(2N)!}.
\end{equation*}
Now, we apply Stirling's formula $\frac{n^n\sqrt{2\pi n}}{n!e^n}\xlongrightarrow{n\to\infty}1$
to approximate, for sufficiently large $N$:
\begin{equation*}
P_{1,N}\approx \frac{N^N\sqrt{2\pi N}(2/e)^N}{(2N)^{2N}\sqrt{4\pi N}/e^{2N}}\\=\frac{1}{\sqrt{2}}\left(\frac{e}{2N}\right)^N.
\end{equation*}
Let $Z_{d,2N}^{p,c}$ be similarly defined to $Z_{d,2N}^p$ but restricting the hyperplanes to be orthogonal to one of the $d$ canonical axes. Then, for any $d\geq1$, we can bound:
\begin{equation}\label{eq:lowboundpc}
P(Z_{d,2N}^{p,c}=N-1)\leq P(Z_{d,2N}^p=N-1).
\end{equation}
By hypothesis, the $d$ random variables defined as $Z_{1,2N}^p$ over the projection of $\mathcal{D}$ on each canonical axis are i.i.d. to $Z_{1,2N}^p$, so we can write: \begin{align*}
P(Z_{d,2N}^{p,c}> N-1)&=\left[1-P(Z_{1,2N}=N-1)\right]^d\\
&=\frac{N!\,2^N}{(2N)!}\approx\left[1-\frac{1}{\sqrt{2}}\left(\frac{e}{2N}\right)^N\right]^d,
\end{align*}
if $N\gg1$. By \eqref{eq:lowboundpc}, the complementary provides the desired lower bound for $P_{d,N}$.
\end{proof}

\begin{proof}[Proof of \cref{cor:snodewithsep}]
Let $\mathbf{a}\in\mathbb{S}^{d-1}$, $\{b_n\}_{n=1}^{N+1}\subset\mathbb{R}$ and $\tau$ be given by \cref{hyp1}. With no loss of generality, we can assume that $\mathbf{a}=\mathbf{e}_1$ and $\tau$ is the identity permutation. The argument that we will use is similar to the one employed in the proof of \cref{th:SimControl-p}, but now the motion must be longitudinal as well as transverse. It also hinges on the fact that, , 
inside the $n$-th strip \begin{equation*}S_n\coloneqq\{\bfx\in\mathbb{R}^d:-b_{n}<\mathbf{a}\cdot\bfx<-b_{n+1}\},\end{equation*} the equation \eqref{eq:node-shallow} simplifies to \eqref{nodestripn}. The simultaneous control of the data points is achieved inductively, in increasing order of the first coordinates, by appropriately defining the field $\mathbf{w}_n$ associated with each hyperplane $H_n$. Both the base case and the inductive step are established by the following two lemmas, which will be proven later.
\begin{lemma}\label{lem1}
Consider two points $\bfx_1,\mathbf{y}_1\in\mathbb{R}^d$ with $\bfx_1\neq\mathbf{y}_1$. For any $T>0$ and $b\in\mathbb{R}$ satisfying $x_1^{(1)}+b>0$ and $y_1^{(1)}+b>0$, there exists a unique $\mathbf{w}\in\mathbb{R}^d$ such that the solution of
\begin{align}\label{cprob1point}
\begin{cases}
\dot\bfx&=\mathbf{w}\sigma(x^{(1)}+b),\\
\bfx(0)&=\bfx_1\in\mathbb{R}^d
\end{cases}
\end{align}
reaches $\bfx(T)=\mathbf{y}_1$.
\end{lemma}
Having controlled $\bfx_1,\dots,\bfx_{n-1}$ to $\mathbf{y}_1,\dots,\mathbf{y}_{n-1}$ in a time horizon $T>0$, and using parameters $\{\mathbf{w}_i\}_{i=1}^{n-1}\subset\mathbb{R}^d$ and $\{b_i\}_{i=1}^{n-1}\subset\mathbb{R}$ such that $b_{n-1}<\cdots<b_1$ and
\begin{equation*}
    x_i^{(1)}+b_i>0\quad\text{and}\quad y_i^{(1)}+b_i>0\quad\text{for }i=1,\dots,N,
\end{equation*}
steering $\bfx_n$ involves overcoming an autonomous drift field
\begin{equation}\label{eq:drift}\mathbf{d}(\bfx)\coloneqq\sum_{i=1}^{n-1}\mathbf{w}_i\sigma(\mathbf{a}\cdot\bfx+b_i)=\sum_{i=1}^{n-1}\mathbf{w}_i\sigma(x^{(1)}+b_i).\end{equation}
The drift field $\mathbf{d}$ becomes more intense as the first coordinate increases, owing to the characteristics of the ReLU function. However, the following lemma shows that the control is possible:
\begin{lemma}\label{lem2}
    Consider $\bfx_n,\mathbf{y}_n\in\mathbb{R}^d$ with $\bfx_n\neq\mathbf{y}_n$. With the above notation, for any $T>0$   there exists a unique $\mathbf{w}_n\in\mathbb{R}^d$ and some $b_n\in\mathbb{R}$ satisfying \begin{equation*}
    b_{n}<b_{n-1},\quad x_n^{(1)}+b_n>0\quad\text{and}\quad y_n^{(1)}+b_n>0
\end{equation*} such that the solution of the Cauchy problem
\begin{align}\label{cprob2point}
\begin{cases}
\dot\bfx&=\mathbf{d}(\bfx)+\mathbf{w}_n\sigma(x^{(1)}+b_n),\\
\bfx(0)&=\bfx_n,
\end{cases}
\end{align}
where $\mathbf{d}$ is given by \eqref{eq:drift}, reaches $\bfx(T)=\mathbf{y}_n$.
\end{lemma}
With \cref{lem1,lem2}, the inductive argument is almost complete. It is left to show that the trajectory of each initial datum $\bfx_n$ will remain in its corresponding strip $S_n$ for all $t\in(0,T)$.  

On one hand, \cref{lem2} guarantees that the trajectory $\bfx(t)$ originating from $\bfx_n$ will reach the endpoint $\mathbf{y}_n$. On the other hand, taking into account that the field is autonomous and invariant along each hyperplane $x^{(1)}=\text{const}$, then $\dot x^{(1)}(t)$ cannot change sign at any time, that is to say, the $x^{(1)}(t)$ does not change its direction. Consequently, the entire trajectory will be contained within the strip bounded by $x=x_n^{(1)}$ and $x=y_n^{(1)}$, which in turn is contained in $S_n$.
\end{proof}

\begin{proof}[Proof of \cref{lem1}]
In the half-space $\{x^{(1)}+b>0\}$, the equation \eqref{cprob1point} is written as
\begin{equation}\label{eq:lem1rest}
\dot\bfx=\mathbf{w}\sigma(x^{(1)}+b)=\mathbf{w} x^{(1)}+\mathbf{w}\, b.
\end{equation}
We can assume that $b=0$, so the solution of \eqref{eq:lem1rest} is
\begin{equation*}
\bfx(t)=\frac{x_1^{(1)}}{w^{(1)}}\mathbf{w}\left(e^{w^{(1)}t}-1\right)+\bfx_1,
\end{equation*}
which can be driven to $\bfx(T)=\mathbf{y}_1$ by taking 
\begin{equation*}
w^{(1)}=\frac{1}{T}\ln\left(\frac{y_1^{(1)}}{x_1^{(1)}}\right)\quad\text{and}\quad w^{(k)}=\frac{y_1^{(k)}-x_1^{(k)}}{y_1^{(1)}-x_1^{(1)}}\,w^{(1)},
\end{equation*}
or $w^{(k)}=y_1^{(k)}-x_1^{(k)}$ if $x_1^{(1)}=y_1^{(1)}$, for $k=2,\dots,d$. Moreover, $\bfx(t)$ stays in the half-space $x^{(1)}>0$ for $t\in[0,T]$ because $x^{(1)}(t)$ is monotone in that interval.
\end{proof}
\begin{proof}[Proof of \cref{lem2}]
First, take any $b_n\in\mathbb{R}$ satisfying \begin{equation*}
    b_{n}<b_{n-1},\quad x_n^{(1)}+b_n>0\quad\text{and}\quad y_n^{(1)}+b_n>0.
\end{equation*}
For simplicity, we rewrite \eqref{cprob1point} as
\begin{equation*}
\dot\bfx=\left(\mathbf{s}_{n-1}+\mathbf{w}_n\right)x^{(1)}+\mathbf{w}_nb_n
+\mathbf{c}_{n-1},
\end{equation*}
where $\mathbf{s}_{n-1}=\sum_{i=1}^{n-1}\mathbf{w}_i$, $\mathbf{c}_{n-1}=\sum_{i=1}^{n-1}\mathbf{w}_ib_i$. If we restrict to the first coordinate, we have:
\begin{align*}
\begin{cases}
\dot x^{(1)}&=\left(s_{n-1}^{(1)}+w_n^{(1)}\right)x^{(1)}+w_n^{(1)}b_n
+c_{n-1}^{(1)},\\[4pt]
x^{(1)}(0)&=x_n^{(1)},
\end{cases}
\end{align*}
which has solution\begin{multline*}x^{(1)}(t)=\frac{w_n^{(1)}b_n+c_{n-1}^{(1)}}{s_{n-1}^{(1)}+w_n^{(1)}}\left[e^{\left(s_{n-1}^{(1)}+w_n^{(1)}\right)t}-1\right]\\+x_n^{(1)}e^{\left(s_{n-1}^{(1)}+w_n^{(1)}\right)t}.\end{multline*}
First, we want to see if there exists $\hat w_n^{(1)}\in\mathbb{R}$ such that $x^{(1)}(T)=y_n^{(1)}$, or, equivalently, if the function \begin{multline}
    f(z)=\frac{z\,b_n+c_{n-1}^{(1)}}{s_{n-1}^{(1)}+z
    }\left[e^{\left(s_{n-1}^{(1)}+z\right)T}-1\right]\\
    +x_n^{(1)}e^{\left(s_{n-1}^{(1)}+z\right)T}-y_n^{(1)}
\end{multline}
has a real root. For that task, we compute:
\begin{align*}
    \lim_{z\to\infty} f(z)&=\operatorname{sign}\left(b_n+x_n^{(1)}\right)\cdot\infty=+\infty,\\
    \lim_{z\to-\infty} f(z)&=-b_n-y_n^{(1)}<0.\\
\end{align*}
On the other hand,
\begin{align*}
    \lim_{z\to{-s_{n-1}^{(1)}}^+} f(z)&=\big(c_{n-1}^{(1)}-b_ns^{(1)}_{n-1}\big)T+x_n^{(1)}-y_n^{(1)}\\&=\lim_{z\to{-s_{n-1}^{(1)}}^-} f(z),
\end{align*}
so $f$ is continuous. Hence, we can assure that there exists $ \hat w_n^{(1)}\in\mathbb{R}$ such that $f\left(\hat w_n^{(1)}\right)=0$. Now, we denote \begin{equation*}s_n^{(1)}=\hat w_n^{(1)}+s_{n-1}^{(1)}\quad\text{and}\quad c_n^{(1)}=\hat w_n^{(1)}b_n+c_{n-1}^{(1)}.\end{equation*} For each component $j\in\{2,\dots,d\}$, 
\begin{align*}
\begin{cases}
\dot x^{(j)}&= c_n^{(1)}\frac{s_{n-1}^{(j)}+w_n^{(j)}}{ s_n^{(1)}}\left[e^{ s_n^{(1)}t}-1\right]\\
&\;+x_n^{(1)}\left(s_{n-1}^{(j)}+w_n^{(j)}\right)e^{s_n^{(1)}t}+w_n^{(j)}b_n+c_{n-1}^{(j)},\\[8pt]
x^{(j)}(0)&=x_n^{(j)},
\end{cases}
\end{align*}
which has solution
\begin{multline*}
x^{(j)}(t)=x_n^{(j)}+\left[w_n^{(j)}b_n+c_{n-1}^{(j)}-c_n^{(1)}\frac{s_{n-1}^{(j)}+w_n^{(j)}}{s_n^{(1)}}\right]t\\
+\left[ c_{n}^{(1)}\frac{s_{n-1}^{(j)}+w_n^{(j)}}{ s_n^{(1)^2}}+x_n^{(1)}\frac{s_{n-1}^{(j)}+w_n^{(j)}}{s_n^{(1)}}\right]e^{s_n^{(1)}t}.
\end{multline*}
Now, we want to find a solution of $g(z)=0$ for \begin{multline*}g(z)=x_n^{(j)}+\left[zb_n+c_{n-1}^{(j)}- c_n^{(1)}\frac{s_{n-1}^{(j)}+z}{s_n^{(1)}}\right]T\\+\left[c_n^{(1)}\frac{s_{n-1}^{(j)}+z}{ s_n^{(1)^2}}+x_n^{(1)}\frac{s_{n-1}^{(j)}+z}{s_n^{(1)}}\right]e^{s_n^{(1)}T}-y_n^{(j)}.\end{multline*}
This is an affine function in $z$, so $g(z)=0$ has a unique solution $\hat w_n^{(j)}\in\mathbb{R}$ if and only if the slope is non-zero.  Suppose that the $b_n$ we have chosen yields a zero slope, i.e.,
\begin{equation*}
\left[\frac{c_n^{(1)}}{s_n^{(1)}}+x_n^{(1)}\right]e^{s_n^{(1)}T}+(b_ns_{n-1}^{(1)}-c_{n-1}^{(1)})T=0.
\end{equation*}
Recalling that
 $f(\hat w_n^{(1)})=0$, we can write this equation as
\begin{equation*}
y_n^{(1)}+\frac{c_n^{(1)}}{s_n^{(1)}}+(b_ns_{n-1}^{(1)}-c_{n-1}^{(1)})T=0,
\end{equation*}
so \begin{equation}\label{eq:eqproof}\frac{c_n^{(1)}}{s_n^{(1)}}=(c_{n-1}^{(1)}-b_ns_{n-1}^{(1)})T-y_n^{(1)}\end{equation} and we can solve for $\hat w_n$ in this equation as \begin{equation*}
\hat w_n^{(1)}=\frac{\left[(c_{n-1}^{(1)}-b_ns_{n-1}^{(1)})T-y_n^{(1)}\right]s_{n-1}^{(1)}-c_{n-1}^{(1)}}{b_n-\left(c_{n-1}^{(1)}-b_ns_{n-1}^{(1)}\right)T-y_n^{(1)}}.
\end{equation*}
Meanwhile, substituting \eqref{eq:eqproof} in $f(\hat w_n^{(1)})=0$, we get:
\begin{equation*}
\left[x_n^{(1)}-y_n^{(1)}\right]e^{s_n^{(1)}T}-T(c_{n-1}^{(1)}+b_ns_{n-1}^{(1)})[e^{s_n^{(1)}T}-1]=0,
\end{equation*}
which can also be solved for $\hat w_n^{(1)}$ as
\begin{equation*}
\hat w_n^{(1)}=\frac{1}{T}\ln\left(\frac{(c_{n-1}^{(1)}-b_ns_{n-1}^{(1)})T}{(c_{n-1}^{(1)}-b_ns_{n-1}^{(1)})T-y_n^{(1)}+x_n^{(1)}}\right)-s_{n-1}^{(1)}.
\end{equation*}
Equalizing both expressions of $\hat w_n^{(1)}$, we obtain an equation with different analytic functions of $b_n$ and not involving $\hat w_n^{(1)}$. Therefore, changing slightly the chosen value of $b_n$, these expressions become different, independently of the corresponding value of $\hat w_n^{(1)}$, so there is a perturbation of $b_n$ that satisfies the statement of the lemma while $g(z)=0$ has a unique solution.

\end{proof}
\subsection{Proof of \cref{th:control-shallow}}
First, we prove \cref{prop:exactsmooth} with an inductive argument of topological nature.
\begin{proof}[Proof of \cref{prop:exactsmooth}]
First, we aim to build a family of $N$ disjoint $C^\infty$ curves contained in $\operatorname{Int}(\Omega)$, each connecting the two points of a corresponding pair $(\bfx_n,\mathbf{y}_n)\in\mathcal{D}$. When $d\geq 2$, any connected open set in $\mathbb{R}^d$ is path-connected. Therefore, we can take a continuous path $C$ that connects any two given points $(\bfx,\mathbf{y})$ inside $U\coloneqq\operatorname{Int}(\Omega)\setminus K$, where $K$ represents any finite union of disjoint curves contained in $\operatorname{Int}(\Omega)$. Moreover, by a well-known approximation argument this path can be chosen to be $C^\infty$. 

Now, we have $N$ disjoint $C^\infty$ curves $\{C_n\}_{n=1}^N$ contained in $\operatorname{Int}(\Omega)$, each connecting a corresponding pair of points $(\bfx_n,\mathbf{y}_n)$. The tangent velocity field of each curve is also $C^\infty$, so the assembled field $\mathbf{V}'$, defined in $\bigcup_{n=1}^N C_n$, is smooth too. Since its domain is compact, it is also Lipschitz-continuous. The required vector field $\mathbf{V}$ is provided by Kirzsbraun's Theorem (see \cite{valentine45}), which ensures the existence of a Lipschitz-continuous map $\mathbf{V}:\mathbb{R}^d\to\mathbb{R}^d$ that extends $\mathbf{V}'$ sharing the same Lipschitz constant.
\end{proof}

The proof of \cref{th:control-shallow} employs the following lemma from \cite[Section~7.2.2]{devore_hanin_petrova_2021}:
\begin{lemma}[Approximation rate for Lipschitz functions]\label{lem:devore}
    Let $K$ be the unit ball in $\operatorname{Lip}(\Omega,\mathbb{R})$, where $\Omega=[-R,R]^d$. We have
    \begin{equation*}
    \frac{1}{[\kappa\log_2\kappa]^{1/d}}\,C_{d,R}\leq E(K,\Sigma_\kappa)_{C(\Omega)}\leq C_{d,R}\,\frac{\log_2\kappa}{\kappa^{1/d}},
    \end{equation*}
    where:
    \begin{itemize}
        \item $\Sigma_\kappa\coloneqq\left\{S_\kappa:\mathbb{R}^d\to\mathbb{R}\right\}$ is the space of shallow neural networks with  $\kappa=(d+2)p$ parameters, $p$ being the number of neurons in the hidden layer;
        \item $E(K,\Sigma_\kappa)_{C(\Omega)}=   \sup_{f\in K}\inf_{S\in\Sigma_\kappa}\|f-S\|_{C(\Omega)}$ measures the capacity of $\Sigma_\kappa$ to approximate any function in $K$.
    \end{itemize}
\end{lemma}
This result was in turn derived from \cite[Proposition 6]{bachbreaking}, where the upper bound is obtained for $K_L$, the ball of radius $L$ in $\operatorname{Lip}(\Omega,\mathbb{R})$: \begin{equation}\label{ratebounds}
   E(K_L,\Sigma_\kappa)_{C(\Omega)}\leq C_{d,R} \,L\;\frac{\log_2\kappa}{\kappa^{1/d}}.
\end{equation}
\begin{proof}[Proof of \cref{th:control-shallow}]
   \Cref{prop:exactsmooth} ensures that we can find a  Lipschitz-continuous field $\mathbf{V}:\mathbb{R}^d\to\mathbb{R}^d$ such that the flow $\Psi_T$ of the ODE \begin{equation}\label{eq1}\dot\bfx=\mathbf{V}(\bfx)\end{equation}
interpolates the dataset. The classical UAT result in \cite{cybenko} guarantees that we can uniformly approximate in $\Omega$ with precision $\varepsilon/\sqrt{d}$ each of its components $V^{(i)}$ using a corresponding shallow neural network $S^i_{\kappa_i}:\Omega\rightarrow\mathbb{R}$, with $\kappa_i\geq1$ for $i=1,\dots,d$. Moreover, \cref{lem:devore} quantifies the dependence of the error with respect to the number of parameters of each $S^i_{\kappa_i}$, ensuring that \begin{equation*}\kappa_1=\kappa_2=\cdots=\kappa_d=(d+2)p\end{equation*}parameters suffice to ensure \eqref{ratebounds} on each component.
    The assembled field \begin{equation*}\mathbf{V}_{NN}=\left(S_{\kappa_1}^1,\dots,S_{\kappa_d}^d\right):\mathbb{R}^d\to\mathbb{R}^d\end{equation*}
   is of the form $\mathbf{V}_{NN}(\mathbf{x})=W\boldsymbol{\sigma}(A\mathbf{x}+\mathbf{b})$ with complexity $\kappa=(d+2)pd$, and satisfies $\sup_{\bfx\in\Omega}|\mathbf{V}(\bfx)-\mathbf{V}_{NN}(\bfx)|<\varepsilon.$ Consider the neural ODE given by \begin{equation}\label{eq2}
       \dot\bfx=\mathbf{V}_{NN}(\bfx).
   \end{equation}
  Let $\mathbf{X}_{V}(t;\bfx_0)$ and $\mathbf{X}_{NN}(t;\bfx_0)$ be the respective trajectories in time $t>0$ that a point $\bfx_0\in\Omega$ will follow under the dynamics provided by \eqref{eq1} and \eqref{eq2}.
   The deviation \begin{equation*}\mathbf{z}(t)=|\mathbf{X}_{V}(t;\bfx_0)-\mathbf{X}_{NN}(t;\bfx_0)|\end{equation*}is bounded as
    \begin{align*}
    \mathbf{z}(t)& \leq \int_{0}^t\left|\mathbf{V}(\mathbf{X}_{V}(s;\bfx_0))-\mathbf{V}_{NN}(\mathbf{X}_{NN}(s;\bfx_0))\right|\mathrm{ds}\\
    &\leq \int_{0}^t\big\{\left|\mathbf{V}(\mathbf{X}_{V}(s;\bfx_0))-\mathbf{V}(\mathbf{X}_{NN}(s;\bfx_0))\right|\\
    &\;+\left|\mathbf{V}(\mathbf{X}_{NN}(s;\bfx_0))-\mathbf{V}_{NN}(\mathbf{X}_{NN}(s;\bfx_0))\right|\big\}\mathrm{ds}\\
    &\leq L_V\int_{0}^t\left|\mathbf{X}_{V}(s;\bfx_0)-\mathbf{X}_{NN}(s;\bfx_0)\right|\mathrm{ds}+\varepsilon\;t\\
&=L_V\int_0^t\mathbf{z}(s)\mathrm{ds}+\varepsilon\;t.
    \end{align*}
    By Grönwall's inequality, it follows that
    \begin{equation*}
    \mathbf{z}(t)\leq \varepsilon\;t\;\exp\left\{L_V\;t\right\}.
    \end{equation*}
On the other hand, in the second line, we could alternatively add and subtract $\mathbf{V}_{NN}(\mathbf{X}_{V}(s;\bfx_0))$.  Then, if we denote by $L_{NN}$ the Lipschitz constant of $\mathbf{V}_{NN}$, we have:
    \begin{align*}
    \mathbf{z}(t)& \leq \int_{0}^t\big\{\left|\mathbf{V}(\mathbf{X}_{V}(s;\bfx_0))-\mathbf{V}_{NN}(\mathbf{X}_{V}(s;\bfx_0))\right|\\
    &\;+\left|\mathbf{V}_{NN}(\mathbf{X}_{V}(s;\bfx_0))-\mathbf{V}_{NN}(\mathbf{X}_{NN}(s;\bfx_0))\right|\big\}\mathrm{ds}\\
    &\leq\varepsilon\;t+ L_{NN}\int_{0}^t\left|\mathbf{X}_{NN}(s;\bfx_0)-\mathbf{X}_{V}(s;\bfx_0)\right|\mathrm{ds}\\&=\varepsilon\;t+L_{NN}\int_0^t\mathbf{z}(s)\mathrm{ds}.
    \end{align*}
    By Grönwall's inequality, it follows that
    \begin{equation*}
    \mathbf{z}(t)\leq \varepsilon\;t\;\exp\left\{L_{NN}\;t\right\}.\end{equation*}
    Taking $\bfx_0=\bfx_n$ and $t=T$, the two bounds for $\mathbf{z}(t)$ give:
    \begin{equation*}
    |\mathbf{y}_n-\Phi_T(\bfx_n)|\leq \varepsilon\; T\;\exp\big\{\min\{L_V,L_{NN}\}\, T\big\}.
    \end{equation*}
    Note that $L_{NN}\leq \|W\|\cdot\|A\|$ because $\boldsymbol{\sigma}$ is $1$-Lipschitz, so the
    approximation rate \eqref{eq:ratethm2} is obtained by direct application of \eqref{ratebounds}.
\end{proof}

    \subsection{Proof of \cref{th:NTEq-control-p}}\label{appendix:NTEq-control}
\begin{figure}[t]  
\centering
    \includegraphics[scale=0.35]{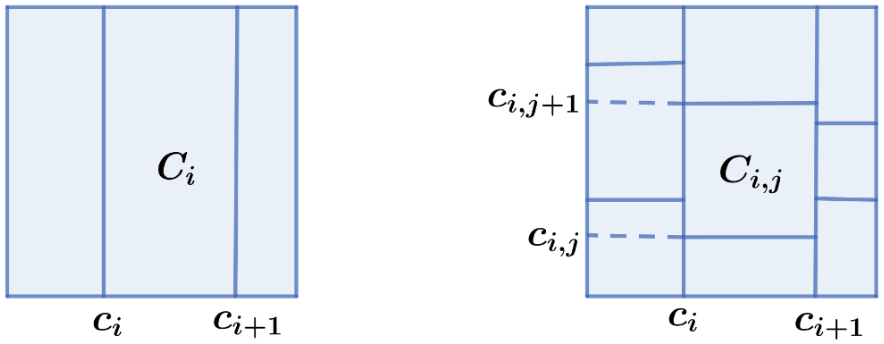}
   \caption{Division of $[0,1]^2$ into rectangles, each containing a mass of $1/n^2$ following the distribution given by $\rho_0$.}
        \label{fig:NTEq-control1}
\end{figure}

\begin{figure}[t]  
\centering
    \includegraphics[scale=0.35]{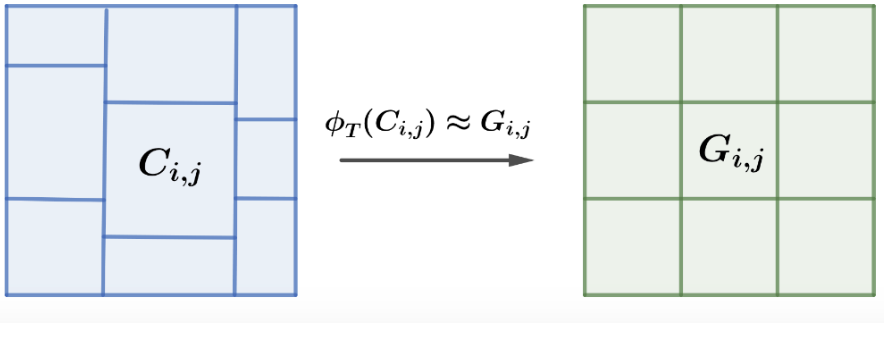}
   \caption{Representation of the ideal transformation of the rectangles $C_I$ into  the corresponding ones $G_I$.}
        \label{fig:NTEq-control1.5}
\end{figure}

\begin{figure}[t]  
\centering
    \includegraphics[width=0.8\linewidth]{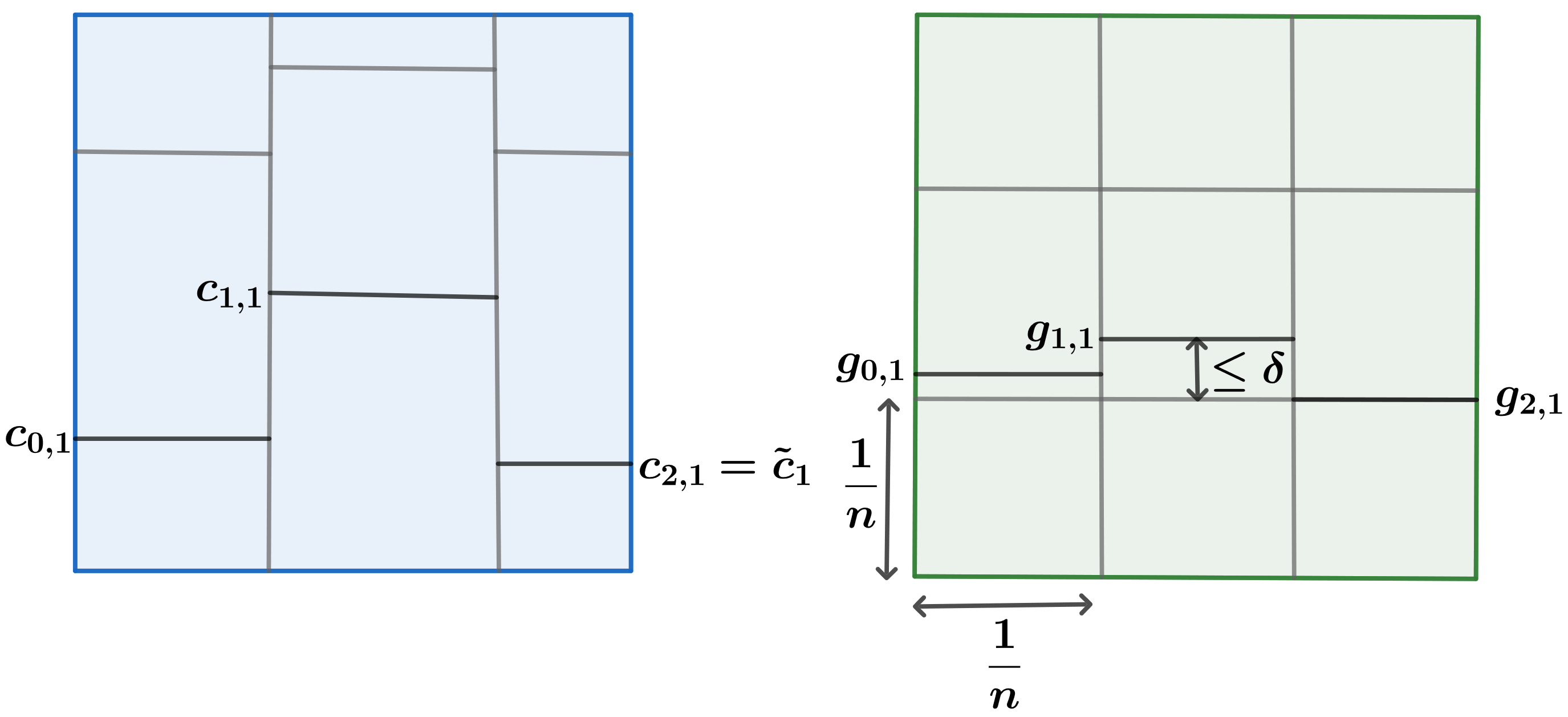}
   \caption{Construction of the partition in rectangles $G_I^\delta$.}
        \label{fig:epsrect}
\end{figure}

       We seek to find $(W,A,\mathbf{b})$ such that the generated vector field moves, compresses and stretches the mass distributed following $\rho_0$, to drive it approximately to the target density, given by $\rho_*$.   This is achieved in four steps, illustrated in \cref{fig:NTEq-control1,fig:NTEq-control1.5,fig:epsrect,fig:NTEq-control2} for the case $d=2$.
    \newline
        \textbf{1. Preparation.} 
         We compress $\operatorname{supp}(\mu_0)$ into $[0,1]^d$. To do this, we aim to   find a control  \begin{equation*}\left(W,A,\mathbf{b}\right)\in L^\infty\left((0,T);\mathbb{R}^{d\times p}\times\mathbb{R}^{p\times d}\times\mathbb{R}^p\right)\end{equation*} such that, in a time $T_1>0$, the flow $\Phi_{T_1}$ of \eqref{eq:node-p} satisfies
        \begin{equation*}\Phi_{T_1}(\operatorname{supp}(\mu_0))\subset [0,1]^d\end{equation*}  
        For $k=1,\dots,d$, we fix the hyperplane $x^{(k)}=0$ and a compressive velocity field by taking $(\mathbf{w},\mathbf{a},b)=(-\mathbf{e}_k,\mathbf{e}_k,0)$. We choose $T_{1,k}>0$ sufficiently large to ensure
\begin{equation*}\Phi_{T_{1,k}}\left(\operatorname{supp}(\mu_0)\cap\{x^{(k)}\geq0\}\right) \subset\{0\leq x^{(k)}\leq1\}.\end{equation*}  
               We repeat the operation with the hyperplanes  $x^{(k)}=1$ for $k=1,\dots,d$,  taking $(\mathbf{w},\mathbf{a},b)=(\mathbf{e}_k,-\mathbf{e}_k,1)$ and $T_{1,k}'>0$ such that \begin{equation*}\Phi_{T_{1,k}'}\circ\Phi_{T_{1,k}}\left(\operatorname{supp}(\mu_0)\cap\{x^{(k)}\leq0\}\right) \subset\{0\leq x^{(k)}\leq1\}.\end{equation*} Both operations are possible in a finite time because $\mu_0$ has compact support. In the end, we will have built piecewise constant controls $(\mathbf{w},\mathbf{a},b)$ that take $2d$ values, such that
                \begin{equation*}\Phi_{T_1}(\operatorname{supp}(\mu_0))\subset [0,1]^d,\quad\text{for } T_1\coloneqq\sum_{k=1}^d\left(T_{1,k}+T_{1,k}'\right).\end{equation*}
                Using $p$ neurons, we can simultaneously apply $p$ controls, because the characteristic curves of \eqref{eq:NTEq} when $(\mathbf{w},\mathbf{a},b)=(\pm\mathbf{e}_k,\pm\mathbf{e}_k,1)$ are parallel to the hyperplanes $\{x^{(l)}=0\}$ for every $l\neq k$. So, the total number of values taken by $(W,A,\mathbf{b})$ is $\lceil 2d/p\rceil$.\newline
                \textbf{2. Partition.}  We aim to divide $[0,1]^d$ into a collection of $n^d$ hyperrectangles, each containing a mass of $1/n^d$, as distributed by $\mu_0$.  The process can be visualized in \cref{fig:NTEq-control1} for $d=2$.  For simplicity, we redefine $\mu_0:=\Phi_{T_1\#}\mu_0$ with density $\rho_0$, now satisfying $\operatorname{supp}(\mu_0)\subset [0,1]^d$. Let $n\geq 1$ and consider the function
                \begin{equation*}t\longmapsto \int_{[0,t)\times [0,1]^{d-1}} d\mu_0=\int_{[0,t)\times [0,1]^{d-1}} \rho_0.\end{equation*}
               This function is continuous, strictly increasing (by absolute continuity), equal to $0$ at $t=0$ and equal to $1$ at $t=1$. Therefore, we can choose $n+1$ numbers
\begin{equation*}c_0=0<c_1<\dots<1=c_n\end{equation*}
such that, for $i_1=0,\dots,n-1,$
\begin{equation*}\int_{[c_{i_1},c_{i_1+1}]\times[0,1]^{d-1}}\rho_0 =\frac{1}{n}.\end{equation*} 
                Similarly, for each $i_1=0,\dots,n-1$ we choose $n+1$ numbers \begin{equation*}c_{i_1,0}=0<c_{i_1,1}<\dots<1=c_{i_1,n}\end{equation*}
   such that, for $i_2=0,\dots,n-1,$\begin{equation*}\int_{[c_{i_1},c_{i_1+1}]\times[c_{i_1,i_2},c_{i_1,i_2+1}]\times[0,1]^{d-2}}\rho_0 =\frac{1}{n^2}.\end{equation*}  
   
                 Repeating this operation recursively for each coordinate, we end up with $n^d$ hyperrectangles \begin{equation*}C_{i_1,\dots,i_d}^0\coloneqq\left[c_{i_1},c_{i_1+1}\right]\times\cdots\times \left[c_{i_1,\dots,i_d},c_{i_1,\dots,i_d+1}\right]\subset\mathbb{R}^d,\end{equation*}
                 with $i_k\in\{0,\dots,n-1\}$ for every $k=1,\dots,d$, such that
\begin{equation*}\int_{C_{i_1,\dots,i_d}^0}\rho_0=\frac{1}{n^d}.\end{equation*}    
                The analogous partition for the uniform measure $\mu_*$ is
\begin{equation*}G_{i_1,\dots,i_d}\coloneqq\left[\frac{i_1}{n},\frac{i_1+1}{n}\right]\times\dots\times\left[\frac{i_d}{n},\frac{i_d+1}{n}\right].\end{equation*}
For the sake of readability, we will denote each multi-index by $I=(i_1,\dots,i_d)\in\{0,\dots,n-1\}^{d}$ and write $C_I^0$ and $G_I$.\newline
                \textbf{3. Control.} We aim to define the controls that expand and compress the mass until the hyperrectangles $C_I^0$ approximate a corresponding collection of $n^d$ hypercubes, each of them containing the same mass $1/n^d$, as distributed by $\mu_*$. Ideally, we would build $\left(W,A,\mathbf{b}\right)$ such that the flow of the ODE \eqref{eq:node-p} satisfied
                \begin{equation}\label{eq:NTEq-control-rect-2}\Phi_T(C_I^0;W,A,\mathbf{b})=G_I
                \end{equation} 
               for each $I\in\{0,\dots,n-1\}^d$. This would be done by transforming each hyperplane $\{x^{(k)}=c_{i_1,...,i_k}\}$ into a target hyperplane $\{x^{(k)}=i_k/n\}$.   
               \newline
                    However, this task is not possible in general, since any $\{x^{(k)}=i_k/n\}$ can be a target for multiple distinct hyperplanes $\{x^{(k)}=c_{i_1,...,i_k}\}$. 
                    We therefore relax the problem to $\delta$-approximate control by considering $\delta$-displacements of the target hyperplanes $\{x^{(k)}=i_k/n\}$. Now we aim to control each $C_I^0$ to a corresponding target  $G_I^\delta$ that is $\delta$-close to $G_I$, for a sufficiently small $\delta>0$. 
                
                    Let us build the new hyperrectangles $G_I^\delta$, a process shown in \cref{fig:epsrect}. For each $k\in\{2,\dots,d\}$,  $(i_1,\dots,i_{k-1})\in\{0,\dots,n-1\}^{k-1}$ and $i_k\in\{0,\dots,n\}$, we define           \begin{equation*}g_{i_1,\dots,i_k}^\delta\coloneqq i_k/n+\delta(c_{i_1,\dots,i_k}-\Tilde{c}_{i_{k}}),\end{equation*}
                where \begin{equation*}\Tilde{c}_{i_{k}}\coloneqq\min \{ c_{i_1',\dots,i_{k-1}',i_k}:(i_1',\dots,i_{k-1}')\in\{0,\dots,n-1\}^{k-1} \}.\end{equation*} Note that $g_{i_1,\dots,i_k}^\delta=1$ whenever $i_k=n$. By construction,
\begin{equation*}g_{i_1,\dots,i_{k-1},i_k}^\delta<g_{i_1',\dots,i_{k-1}',i_k}^\delta \iff c_{i_1,\dots,i_{k-1},i_k}<c_{i_1',\dots,i_{k-1}',i_k}\end{equation*}
and
\begin{equation*}g_{i_1,\dots,i_{k-1},i_k}^\delta=g_{i_1',\dots,i_{k-1}',i_k}^\delta \iff c_{i_1,\dots,i_{k-1},i_k}=c_{i_1',\dots,i_{k-1}',i_k}.\end{equation*}
   By recursion, we define a new partition of $[0,1]^d$ into a collection of rectangles $G_I^\delta$ with $I\in\{0,\dots,n-1\}^d$, where \begin{equation*}G_I^\delta\subset G_I+\{0\}\times[-\delta,\delta]^{d-1}.\end{equation*}
Moreover, this partition mimics the structure of the partition defined for $\mu_0$, in the sense that there is the same number of distinct target hyperplanes as initial ones to be controlled.  To sum up, if we take $\delta<1/n$, we end up with:
                \begin{align}
                    \label{eq:hyp-subdivision}
                    \left\{C_I^0 : I\in\{0,\dots,n-1\}^d\right\} \text{ s.t. } &\int_{C_I^0}d\mu_0=\frac{1}{n^d},\\\notag
                    \left\{G_I^\delta : I\in\{0,\dots,n-1\}^d\right\} \text{ s.t. }&\int_{G^\delta_I}d\mu_*\leq \frac{3^d}{n^d},\\\notag 
                    \text{ and }&\operatorname{diam}(G^\delta_I)\leq 3\sqrt{d}\frac{1}{n}.
                \end{align}          
            It is left to map $C_I^0$ to $G_I^\delta$ for each $I$. This is based on the following lemma, whose proof we postpone for readability. 
                \begin{lemma}
                    \label{lemma:control-subdiv}
                    Let $d\geq2$, $\mu_0\in\mathcal{P}_{ac}^c(\mathbb{R}^d)$ with density $\rho_0$, $\rho_*$ the uniform density in $[0,1]^d$, and $T>0$ be fixed. Let $n\geq1$ and consider a family of hyperrectangles such as \eqref{eq:hyp-subdivision}.  For any $p_1,\dots,p_d\geq1$, there exists a piecewise constant control \begin{equation*}(W,A,\mathbf{b})\in L^\infty((0,T);\mathbb{R}^{p\times d}\times\mathbb{R}^{d\times p}\times \mathbb{R}^p)\end{equation*}with $p=p_1+\dots+p_d$ such that, for each $I=(i_1,\dots,i_d)\in\{0,\dots,n-1\}^d,$ the flow $\Phi_T$ generated by \eqref{eq:node-p} satisfies
\begin{equation*}\Phi_T(C^0_I;W,A,\mathbf{b})=G^\delta_I\end{equation*}
Furthermore, the number of discontinuities of $(W,A,\mathbf{b})$ is 
\begin{equation*}L=\max\{ \lceil n/p_1\rceil ,\dots, \lceil n^d/p_d\rceil\}-1.\end{equation*}

                \end{lemma}
                
                \textbf{4. Estimates.} We compute the $W_q$-distance between both measures to verify the approximate control.     Let $\Phi_T$ be the flow given by \cref{lemma:control-subdiv}, satisfying
                                   \begin{equation*}\Phi_T(C_I)=G_I^\delta,\; \text{ for } I=(i_1,\dots,i_d)\in\{0,\dots,n-1\}^d.\end{equation*}
                Let us quantify, in the Wasserstein-$q$ distance, the proximity of $\Phi_{T\#}\mu_0$ to $\mu_*$. We have
                \begin{multline}
                \label{eq:W1} W_q(\mu(T),\mu_*)=W_q(\Phi_{T\#}\mu_0,\mu_*)
                    \\\leq\sum_{I\in\{0,\dots,n-1\}^d}W_q(\Phi_{T\#}\mu_0|_{G_I^\delta},\mu_*|_{G_I^\delta}),
                \end{multline}
                see \cite{VILLANI_2016} for the inequality. For each $I\in\{0,\dots,n-1\}^d$, let $\gamma_I:\mathbb{R}^d\rightarrow\mathbb{R}^d$ be the measurable function that satisfies
                \begin{equation*}\gamma_{I\#}(\Phi_{T\#}\mu_0|_{G_I^\delta})=\mu_*|_{G_I^\delta},\end{equation*}
             attaining the minimum in Monge formulation \eqref{eq:monge} for $W_q$. In particular, $\gamma_{I}$ only redistributes the mass inside $G_I^\delta$, so
                \begin{align*}
                    \int_{\mathbb{R}^d}|x-\gamma_{I}(x)|^qd\mu_*|_{G_{I}^\delta}&=\int_{G_{I}^\delta}|x-\gamma_{I}(x)|^q d\mu_*
                    \\
                    &\leq \operatorname{diam}(G_{I}^\delta)^q\int_{G_{I}^\delta}d\mu_*
                    \leq \frac{3^{q+d}d^{q/2}}{n^{q+d}}.
                \end{align*}
                Plugging this bound into \eqref{eq:W1}, we conclude that
               \begin{align*}W_q(\mu(T),\mu_*)&\leq 3^{1+d/q}\sqrt{d}\,n^d\left(\frac{1}{n^{q+d}}\right)^{1/q}\\
               &=3^{1+d/q}\sqrt{d}\, n^{-(1+d/q-d)}.\end{align*}
               It follows that \begin{equation*}W_q(\mu(T),\mu_*)\xlongrightarrow{n\to\infty} 0\;\iff\;q<\frac{d}{d-1},\end{equation*} and, in that case, $W_q(\mu(T),\mu_*)<\varepsilon$ is obtained for \begin{equation*}
n>\left(\frac{3^{1+d/q}\sqrt{d}}{\varepsilon}\right)^{\frac{1}{1+d/q-d}},
                \end{equation*}
                hence proving \cref{th:NTEq-control-p}.
\begin{proof}[Proof of \cref{lemma:control-subdiv}]

We control the rectangles by mapping each hyperplane to its corresponding target hyperplane. This strategy is illustrated in \cref{fig:NTEq-control2}.

The proof is divided into three steps. First, we achieve simultaneous control for any $p$ hyperplanes orthogonal to a fixed direction. Second, for $N \geq p$ hyperplanes, we iteratively control $p$-subsets. Third, we show that this approach is applicable to all $d$ canonical directions simultaneously. \begin{figure}
    \includegraphics[width=0.8\linewidth]{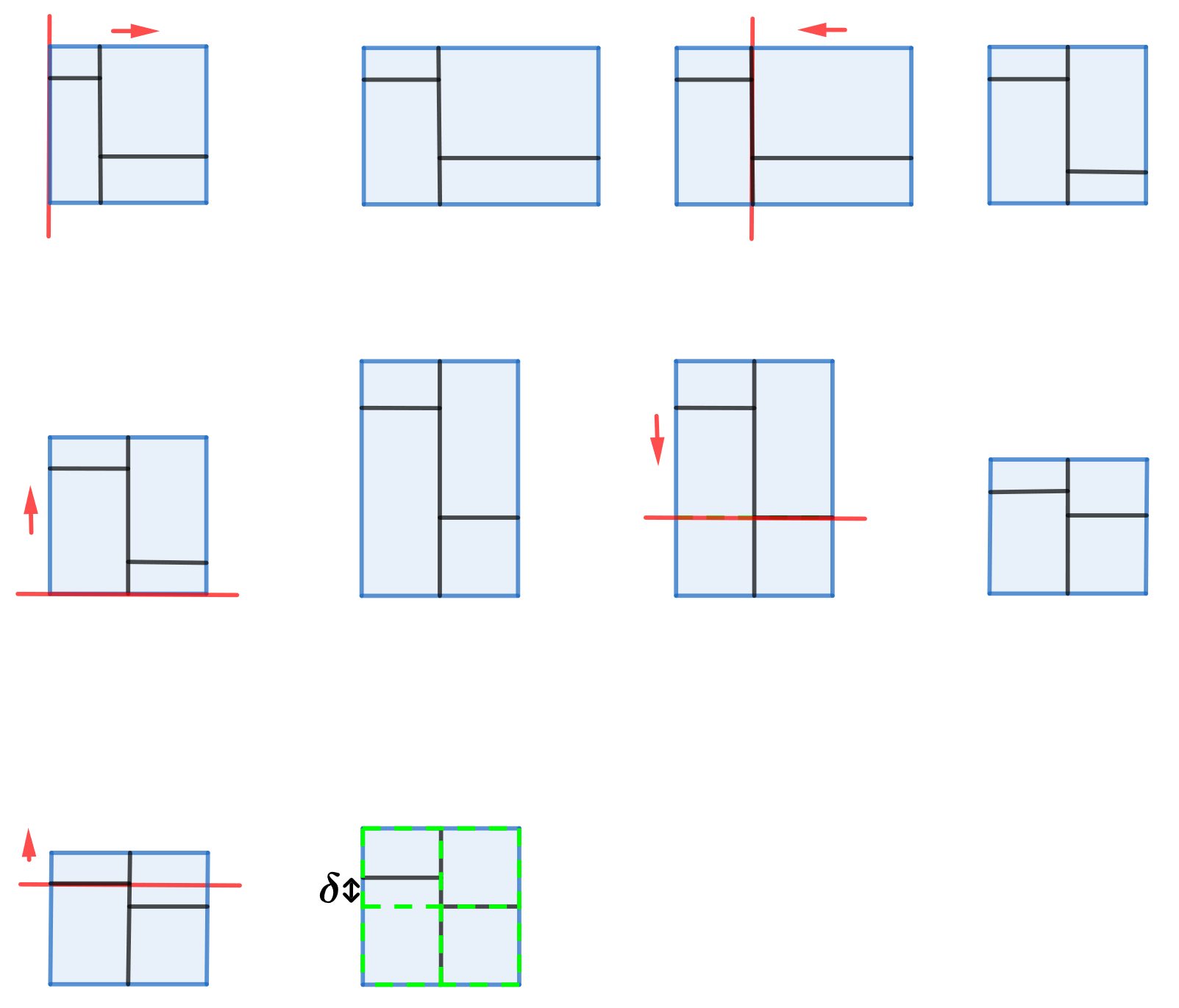}
     \caption{From left to right, by rows: We transform the rectangles $C_I$ to approximate the corresponding ones $G_I^\delta$ (in green)}
          \label{fig:NTEq-control2}
\end{figure}\newline
\textbf{Step 1.} Let $k\in\{1,\dots,d\}$, $p \geq1$ and 
            \begin{align*}
                -\infty<c_1<\cdots<c_p<\infty,\\
                -\infty<g_1<\cdots<g_p<\infty.
            \end{align*}
            We aim to find constant controls $\{(\mathbf{w}_i,\mathbf{a}_i,b_i)\}_{i=1}^p\subset\mathbb{R}^d\times\mathbb{R}^d\times\mathbb{R}$ such that the flow generated by \eqref{eq:node-p} satisfies          
            \begin{align}
            \Phi_T(\{x^{(k)}=c_i\})=\{x^{(k)}=g_i\}, \quad i=1,\dots,p.
            \end{align}
    In particular, we will take $\mathbf{a}_i=\mathbf{e}_k$ for all $i$, and $b_i$ such that $-b_i<\min\{c_i,g_i\}$.
          Since the field \begin{equation*}\sum_{i=1}^p \mathbf{w}_i\sigma(\mathbf{a}_i\cdot\bfx+b_i)=\sum_{i=1}^p\mathbf{w}_i\sigma(x^{(k)}+b_i)\end{equation*} 
          only depends on the $x^{(k)}$-coordinate, it is projectable onto the $x^{(k)}$-axis, i.e., the forward evolution of a hyperplane orthogonal to any coordinate-axis is a hyperplane orthogonal to the same coordinate-axis, for any time. Therefore, we can identify each $\{x^{(k)}=c_i\}$ and $\{x^{(k)}=g_i\}$ with the point $c_i\in\mathbb{R}$ or $g_i\in\mathbb{R}$ and study their evolution in the real line, so the problem becomes one-dimensional. 
          Thus, we identify $x^{(k)}\equiv x$ and fix  $\mathbf{w}_i=w_i \mathbf{e}_k$ with $w_i\in\mathrm{R}$, so we aim to find $(w_i)_{i=1}^p\subset\mathbb{R}$ and $(b_i)_{i=1}^p\subset\mathrm{R}$ such that $\Phi_T(c_i)=g_i$ for $i=1,\dots,p$.
          We proceed by induction on $p$.\newline
         First, we consider $p=1$ and let $c_1,g_1\in\mathbb{R}$.  Take any $-b<\min\{c_1,g_1\}$, so $\{c_{1},g_{1}\}\subset\{x+b>0\}$,  
                    and $w=\operatorname{sign}(g_{1}-c_{1}),$ pointing from $c_1$ to $g_1$.                   Then,
                        \begin{equation*}\Phi_{\Tilde T}\big(c_1\big)=g_1,\quad \text{for }\Tilde T=\left|\log\left(\frac{g_{1}-b}{c_{1}-b}\right)\right|.\end{equation*}
                 Using a time rescaling argument identical to the one in the proof of \ref{th:SimControl-p}, this control can be done in any time $T>0$. Also note that $\{x+b\leq0\}$ is fixed by the flow $\Phi_t$.\newline
        In the inductive step, we assume that the statement is true for some $p\geq1$, and consider
            \begin{align*}
                -\infty<c_1<\cdots<c_{p+1}<\infty,\\
                -\infty<g_1<\cdots<g_{p+1}<\infty.
            \end{align*}
            Let $\{(w_i,b_i)\}_{i=1}^p\subset \mathbb{R}\times\mathbb{R}$ with $-b_i<\min\{c_i,g_i\}$ for all $i$, such that the flow of the one-dimensional problem satisfies
            \begin{equation*}\Phi_T(c_i)=g_i,\quad i=1,\dots,p.\end{equation*}
           We want to add a new pair $(w_{p+1},b_{p+1})\in\mathbb{R}\times\mathbb{R}$ such that 
            \begin{equation*}\Phi_T(c_i)=g_i,\quad i=1,\dots,p+1.\end{equation*}       
                \emph{Case 1: $c_{p}< g_{p}$.}
                    Take $b_{p+1}=-g_{p}$,                    so \eqref{eq:node-p} becomes
                    \begin{equation*}\dot x= \sum_{i=1}^p w_i(x+b_i)\mathbbm{1}_{c_i<x}(x)+w_{p+1}(x-g_{p})\mathbbm{1}_{g_{p}<x}(x).\end{equation*}
                    The added velocity $w_{p+1}(x-g_{p})$ only acts on the half-space $\{x \geq g_p\}$, so the points $\{c_i\}_{i=1}^p$ are only subject to the drift field \begin{equation*}d(x)\coloneqq\sum_{i=1}^p w_i(x+b_i)\mathbbm{1}_{c_i<x}(x).\end{equation*}
                    Therefore, when adding $(w_{p+1},b_{p+1})$, we still have
                    \begin{equation*}\Phi_T(c_i)=g_i,\quad i=1,\dots,p.\end{equation*}
                   If $c_{p+1}\leq g_{p}$, there is $0\leq s<T$ such that $\Phi_s(c_{p+1})=g_{p}$. Otherwise, if $c_{p+1}> g_p$, consider $s=0$. Note that $s$ only depends on the $p$ first neurons, on $c_{p+1}$ and on $g_p$, and it is thus independent of $(w_{p+1},b_{p+1})$. Therefore, $c_{p+1}$ is only subject to $d(x)$ for $t\in(0,s)$, and to $d(x)+w_{p+1}(x-g_{p})$ for $t\in(s,T)$. More precisely: 
                    \begin{equation*}\frac{d}{dt}\Phi_t(c_{p+1})=d(\Phi_t(c_{p+1}))\\
                    + w_{p+1}(\Phi_t(c_{p+1})-g_{p})\mathbbm{1}_{s\leq t<T}(t).\end{equation*}
                    Then, with a similar computation to that of \cref{lem2}:
                    \begin{align*}
                    \Phi_T(c_{p+1})=&\Phi_{T-s}\circ\Phi_{s}(c_{p+1})\\
                        =&\Big(g_p+\frac{\sum_{i=1}^p w_ib_i -w_{p+1} g_{p}}{\sum_{i=1}^{p}w_i+w_{p+1}}\Big)e^{(T-s)\sum_{i=1}^{p+1}w_i}\\
                        &-\frac{\sum_{i=1}^p w_ib_i -w_{p+1} g_{p}}{\sum_{i=1}^{p} w_i+w_{p+1}}.
                    \end{align*}
The expression converges to $g_p$ when $w_{p+1}\to-\infty$, and diverges to $\infty$ when $w_{p+1}\to\infty$.  Therefore, by continuity, and since $g_{p}<g_{p+1}$, there exists $w_{p+1}\in\mathbb{R}$ such that $\Phi_T(c_{p+1})=g_{p+1}$.\newline
                    \emph{Case 2: $g_{p}\leq c_{p}< g_{p+1}$.}
                    Take $b_{p+1}=-c_p$, so   \eqref{eq:node-p} becomes 
                   \begin{equation*}\dot x= \sum_{i=1}^p w_i(x+b_i)\mathbbm{1}_{c_i<x}(x)+w_{p+1}(x-c_{p})\mathbbm{1}_{c_{p}<x}(x),\end{equation*}
                   Again, the points $\{c_i\}_{i=1}^p$ are only subject to $d(x)$ because the added velocity only acts on  $\{x \geq c_p\}$, so \begin{equation*}\Phi_T(c_i)=g_i,\quad i=1,\dots,p.\end{equation*}
                    The point $c_{p+1}$ is subject to the total velocity at $t=0$. A similar computation to Case 1 leads to                \begin{align*}\Phi_T(c_{p+1})=&\left(c_p+\frac{\sum_{i=1}^p w_ib_i -w_{p+1}c_{p}}{\sum_{i=1}^p w_i+w_{p+1}}\right)e^{T\sum_{i=1}^{p+1} w_i}\\
                    &-\frac{\sum_{i=1}^p w_ib_i -w_{p+1} c_{p}}{\sum_{i=1}^p w_i +w_{p+1}}.\end{align*}
                   $\Phi_T(c_{p+1})$ tends to $c_p$ when $w_{p+1}\to-\infty$, and diverges to $\infty$ when $w_{p+1}\to\infty$. By continuity, and since $c_p <g_{p+1}$, there exists $w_{p+1}\in\mathbb{R}$ such that $\Phi_T(c_{p+1})=g_{p+1}$.\newline
                    \emph{Case 3: $g_{p+1}\leq c_{p}$.}  
                  Take $b_{p+1}=-c_p$, so \eqref{eq:node-p} becomes
                    \begin{equation*}\dot x= \sum_{i=1}^p w_i(x+b_i)\mathbbm{1}_{c_i<x}(x)+w_{p+1}(x-c_{p})\mathbbm{1}_{c_{p}<x},\end{equation*}
                    and
    \begin{equation*}\Phi_T(c_i)=g_i,\quad i=1,\dots,p.\end{equation*}
                    Let $0<s\leq T$ be the first time such that
                    $\Phi_s(c_{p+1})=c_{p}$.
                    An analogous computation to the previous cases gives \begin{align*}
                        \Phi_s(c_{p+1})=&\left(c_{p+1}+\frac{\sum_{i=1}^p w_ib_i -w_{p+1}c_p}{\sum_{i=1}^{p+1} w_i }\right)e^{s\sum_{i=1}^{p+1} w_i}\\
                        &-\frac{\sum_{i=1}^p w_ib_i -w_{p+1}c_{p}}{\sum_{i=1}^{p+1} w_i},
                    \end{align*}
                    so
                    \begin{equation*}s=\frac{1}{\sum\limits_{i=1}^{p+1} w_i}\log\left(\frac{\sum\limits_{i=1}^p w_i(c_p-b_i) +2w_{p+1} c_{p}}{\sum\limits_{i=1}^p w_i(c_{p+1}+b_i) +w_{p+1} (c_{p+1}-c_{p})}\right).\end{equation*} 
                    By varying $w_{p+1}$ in $(-\sum_{i=1}^p w_i,+\infty)$, we can ensure that $s$ can take any value in $(0,T)$. Then,
                    \begin{align*}
                        \Phi_T(c_{p+1})=&\Phi_{T-s}\circ\Phi_{s}(c_{p+1})\\
                        =&\left(c_p+\frac{\sum_{i=1}^p w_ib_i }{\sum_{i=1}^p w_i }\right)e^{(T-s)\sum_{i=1}^p w_i }\\
                        &-\frac{\sum_{i=1}^p w_ib_i}{\sum_{i=1}^p w_i}.
                    \end{align*}                    
                   For $s=T$, the expression is equal to $c_p$, and for $s=0$ it equals $\Phi_T(c_p)=g_p$. By continuity, and also because $g_p<g_{p+1}\leq c_p$, an argument like in case 1 ensures the existence of $w_{p+1}\in\mathbb{R}$ such that $\Phi_T(c_{p+1})=g_{p+1}.$\newline
    \textbf{Step 2.} Let $k\in\{1,\dots,d\}$, $p\geq1,$ and $N\geq p$ and 
    \begin{equation*}
        -\infty<c_1<\dots<c_N<\infty,\quad
         -\infty<g_1<\dots<g_N<\infty.
    \end{equation*}
    We will show that there exist piecewise constant controls $\left(w_i,a_i,b_i\right)_{i=1}^p$ such that the flow of \eqref{eq:node-p} satisfies
    \begin{equation*}
        \Phi_T(c_i)=g_i,\quad i=1,\dots,N,
    \end{equation*}
and the number of discontinuities is $\lceil N/p\rceil-1$. We use a similar argument to the one in the proof of \cref{th:SimControl-p}. We divide $\{c_i\}_{i=1}^N$ and $\{g_i\}_{i=1}^N$ into subsets of $p$ points \begin{equation*}C_j\coloneqq\{c_{(j-1)\cdot p+1},\dots,c_{j\cdot p}\},\;G_j\coloneqq\{g_{(j-1)\cdot p+1},\dots,g_{j\cdot p}\},\end{equation*} 
for $j=1,\dots,\lceil N/p\rceil-1$, and $C_{\lceil N/p\rceil}$, $G_{\lceil N/p\rceil}$ with the remaining $N-p\lfloor N/p\rfloor$ points. 

The piecewise constant controls are obtained by induction on $j$. In each iteration, we apply step 1 to define $p$ constant controls $(w_i^j,a_i^j,b_i^j)_{i=1}^p$ that map the $p$ points of $C_j$ to the corresponding ones in $G_j$ in time $\frac{T}{\lceil N/p\rceil}$. Note that the initialization of induction in step 1 ensures that the previously controlled subsets $C_1,\dots,C_{j-1}$ can remain fixed during the subsequent iterations, which trivializes the induction. Finally, we have the piecewise constant controls
    \begin{equation*}(w_i,a_i,b_i)_{i=1}^p=\sum_{j=1}^{\lceil N/p\rceil} (w_i^j,a_i^j,b_i^j)_{i=1}^p\mathbbm{1}_{\left(\frac{(j-1)T}{\lceil N/p\rceil},\frac{jT}{\lceil N/p\rceil}\right)}(t),\end{equation*}
    which achieve the desired objective.\newline
    \textbf{Step 3.}  For every $k=1,\dots,d$, let $p_k,N_k\geq1$, and
    \begin{equation*}
        -\infty<c_1^k<\cdots<c_{N_k}^k<\infty,\quad -\infty<g_1^k<\cdots<g_{N_k}^k<\infty.       
    \end{equation*}
    For each fixed direction $k\in\{1,\dots,d\}$, step 2 is used to build piecewise constant controls \begin{equation*}(\mathbf{w}_{j,k},\mathbf{a}_{j,k},b_{j,k})_{j=1}^{p_k}=(w_{j,k}\,\mathbf{e}_k,\mathbf{e}_k,b_{j,k})_{j=1}^{p_k}\end{equation*} with $\lceil N_k/p_k\rceil-1$ discontinuities, such that  \begin{equation*}\Phi_T(\{x^{(k)}=c_i^k\})=\{x^{(k)}=g_i^k\},\quad i=1,\dots,N.\end{equation*} Moreover, only the $k$-th coordinate is varying on each flow, as argued in steps 1 and 2 when we simplified the problem to one dimension. We define each tern of the assembled control 
$(\mathbf{w}_{j},\mathbf{a}_{j},b_{j})_{j=1}^{p}$, with $p=\sum_{k=1}^d p_k$, by 
\begin{equation*}
(\mathbf{w}_{j},\mathbf{a}_{j},b_{j})=\big(\mathbf{w}_{j-\sum_{i=1}^{k-1}p_{i},p_k},\mathbf{a}_{j-\sum_{i=1}^{k-1}p_{i},p_k},b_{j-\sum_{i=1}^{k-1}p_{i},p_k}\big),
\end{equation*}
for $\sum_{i=1}^{k-1}p_{i}+1\leq j\leq \sum_{i=1}^{k}p_{i}$. Therefore,  
    the resulting neural ODE \eqref{eq:node-p} on each coordinate writes
    \begin{equation*}
     \dot x^{(k)}=w_{1,k}^{(k)}(x^{(k)}+b_{1,k})+\dots+w_{p_k,k}^{(k)}(x^{(k)}+b_{p_k,k})
     \end{equation*}
    All the equations of the system are independent, so each movement does not interfere with the other $d-1$ movements. Therefore, the corresponding flow of \eqref{eq:node-p} satisfies
    \begin{equation*}
        \Phi_T(\{x^{(k)}=c_i^k\})=\{x^{(k)}=g_i^k\}
    \end{equation*}
for $k=1,\dots,d$ and $i=1,\dots,N_k$, and moreover, the number of discontinuities in the controls is \begin{align*}L&=\max\{ \lceil N_1/p_1\rceil-1,\dots, \lceil N_d/p_d\rceil -1\}\\
&=\max_{k=1,\dots,d}\lceil N_k/p_k\rceil-1.\end{align*}  Recalling that $N_k=n^k$ for $k=1,\dots,d$ (by construction of the rectangles $C_I^0$ and $G_I^\delta$), it follows the desired result.
    
\end{proof}

\section{Declaration of competing interest}
The authors declare that they have no known competing financial interests or personal relationships that could have appeared to
influence the work reported in this paper.

\section{Acknowledgments}

This paper was supported by the Madrid Government (Comunidad de Madrid – Spain) under the multiannual Agreement with UAM in the line for the Excellence of the University Research Staff in the context of the V PRICIT (Regional Programme of Research and Technological Innovation). A. Álvarez-López has been funded by a contract FPU21/05673 from the Spanish Ministry of Universities. A. Hadj Slimane has been funded by École Normale Supérieur Paris-Saclay and Université Paris-Saclay. E. Zuazua has been funded by the Alexander von Humboldt-Professorship program, ModConFlex Marie Curie Action, HORIZON-MSCA-2021-DN-01, COST Action MAT-DYN-NET, Transregio 154 Project “Mathematical Modelling, Simulation and Optimization Using the Example of Gas Networks” of the DFG, grants PID2020-112617GB-C22 and TED2021-131390B-I00 of MICINN (Spain).

\bibliographystyle{elsarticle-harv}
\bibliography{biblio.bib}

\end{document}